\documentclass[A4j,11pt]{article}
\usepackage{amsfonts,amsmath,amscd}
\usepackage{graphics}
\usepackage{amssymb}
\usepackage[all]{xy}
\begin{document}

\title{{\bf On the indices of minimal orbits\\
of Hermann actions}}
\author{{\bf Naoyuki Koike}}
\date{}
\maketitle

\begin{abstract}
We give a formula to determine the indices of special (non-totally geodesic) 
minimal orbits of Hermann actions.  
Also, we give examples of such minimal orbits of Hermann actions and 
calculate their indices by using the formula.  
\end{abstract}

\vspace{0.5truecm}






\section{Introduction}
In 1987, Y. Ohnita [O] gave a formula to calculate the indices (and nullities) 
of totally geodesic submanifolds in a symmetric space $N$ of compact type 
and showed that the indices of all Helgason spheres in every simply connected 
irreducible compact symmetric space are equal to zero, 
that is, they are stable.  
In 1993, O. Ikawa [I1] investigated the Jacobi operator of equivariant minimal 
homogeneous submanifold in a Riemannian homogeneous space.  
In 1995, by using Ohnita's index formula, M. S. Tanaka [Ta] determined 
the stability of all polars and meridians in every simply connected 
irreducible compact symmetric space.  
Note that polars and meridians are totally geodesic.  
In 2008, by the index formula, T. Kimura [Ki] determined the stability 
of all totally geodesic singular orbits of all cohomogeneity one actions on 
every simply connected irreducible compact symmetric space.  
In 2009, by using this index formula, T. Kimura and M. S. Tanaka [KT] 
determined the stability of all maximal totally geodesic submanifolds in 
every simply connected irreducible compact symmetric space of rank two.  
Let $N=G/K$ be a symmetrc space of compact type 
equipped with the $G$-invariant metric induced from the Killing form of the Lie 
algebra of $G$.  In this paper, we treat only a symmetric space of compact type equipped 
with such a $G$-invariant metric.  Let $H$ be a symmetric 
subgroup of $G$ (i.e., $({\rm Fix}\,\tau)_0\subset H\subset{\rm Fix}\,\tau$ 
for some involution $\tau$ of $G$), where 
${\rm Fix}\,\tau$ is the fixed point group of $\tau$ and 
$({\rm Fix}\,\tau)_0$ is the identity component of ${\rm Fix}\,\tau$.  
The natural action of $H$ on $N$ is called a {\it Hermann action} 
(see [HPTT], [Kol]).  
Let $\theta$ be an involution of $G$ with $({\rm Fix}\,\theta)_0\subset K
\subset{\rm Fix}\,\theta$.  According to [Co], in the case where $G$ is simple, we may 
assume that $\theta\circ\tau=\tau\circ\theta$ by replacing $H$ to a suitable conjugate 
group of $H$ if necessary except for the following three Hermann action:

\vspace{0.2truecm}

(i) $Sp(p+q)\curvearrowright SU(2p+2q)/S(U(2p-1)\times U(2q+1))
\quad(p\geq q+2)$, 

(ii) $U(p+q+1)\curvearrowright Spin(2p+2q+2)/Spin(2p+1)\times_{{\bf Z}_2}
Spin(2q+1)\quad(p\geq q+1)$, 

(iii) $Spin(3)\times_{{\bf Z}_2}Spin(5)\curvearrowright 
Spin(8)/\omega(Spin(3)\times_{{\bf Z}_2}Spin(5))$, 

\vspace{0.2truecm}

\noindent
where $\omega$ is the triality automorphism of $Spin(8)$.  
Here we note that we remove transitive Hermann actions.  

\vspace{0.5truecm}

\noindent
{\bf Assumption.} 
In the sequel, we assume that 
$\theta\circ\tau=\tau\circ\theta$.  

\vspace{0.5truecm}

Let $\mathfrak g,\mathfrak k$ and $\mathfrak h$ be the Lie algebras of $G,K$ 
and $H$, respectively.  Denote the involutions of $\mathfrak g$ induced form 
$\theta$ and $\tau$ by the same symbols $\theta$ and $\tau$, respectively.  
Set $\mathfrak p:={\rm Ker}(\theta+{\rm id})$ and $\mathfrak q:=
{\rm Ker}(\tau+{\rm id})$.  The vector space $\mathfrak p$ is identified with 
$T_{eK}(G/K)$, where $e$ is the identity element of $G$.  
Take a maximal abelian subspace $\mathfrak b$ of 
$\mathfrak p\cap\mathfrak q$.  
For each $\beta\in\mathfrak b^{\ast}$, we set 
$\mathfrak p_{\beta}:=\{X\in\mathfrak p\,\vert\,{\rm ad}(b)^2(X)=-\beta(b)^2X
\,\,(\forall\,b\in\mathfrak b)\}$ and 
$\triangle':=\{\beta\in\mathfrak b^{\ast}\setminus\{0\}\,\vert\,
\mathfrak p_{\beta}\not=\{0\}\}$.  This set $\triangle'$ is a root system.  
Note that we call $\triangle'$ a root system because $\beta$'s ($\beta\in\triangle'$) give 
a root system in the vector subspace spanned by them (in the sense of [He]) even if 
they do not span $\mathfrak b^{\ast}$.  
Let $\Pi'=\{\beta_1,\cdots,\beta_r\}$ be the simple root system of 
the positive root system $\triangle'_+$ of $\triangle'$ under a lexicographic 
ordering of $\mathfrak b^{\ast}$.  Set 
${\triangle'}^V_+:=\{\beta\in\triangle'_+\,\vert\,\mathfrak p_{\beta}\cap
\mathfrak q\not=\{0\}\}$ and ${\triangle'}^H_+:=\{\beta\in\triangle'_+\,\vert\,
\mathfrak p_{\beta}\cap\mathfrak h\not=\{0\}\}$.  
Define a subset $\widetilde C$ of $\mathfrak b$ by 
$$\begin{array}{l}
\displaystyle{\widetilde C:=\{b\in\mathfrak b\,\vert\,0<\beta(b)<\pi\,
(\forall\,\beta\in{\triangle'}^V_+),\,\,-\frac{\pi}{2}<\beta(b)<
\frac{\pi}{2}\,(\forall\,\beta\in{\triangle'}^H_+)\}.}
\end{array}$$
The closure $\overline{\widetilde C}$ of $\widetilde C$ is a simplicial 
complex.  
Set $C:={\rm Exp}(\widetilde C)$, where ${\rm Exp}$ is the exponential map 
of $G/K$ at $eK$.  
Each principal $H$-orbit passes through only one point of $C$ and 
each singular $H$-orbit passes through only one point of 
${\rm Exp}(\partial\widetilde C)$.  For each simplex $\sigma$ of 
$\overline{\widetilde C}$, only one minimal $H$-orbit through 
${\rm Exp}(\sigma)$ exists.  See proofs of Theorems A and B in [Koi2] 
(also [I2]) about this fact.  
Also, it is konown that only one minimal $H$-orbit through 
${\rm Exp}(\sigma)$ is unstable if $\sigma$ is not a vertex 
(see the proof of Theorem 2.24 in [I2]).  
Denote by $D(H)$ the set of all equivalence classes of (finite dimensional) 
irreducible complex representations of $H$ and $\rho_{G/H}:H\to GL(\mathfrak q)$ the 
isotropy representation of $G/H$, that is, $\rho_{G/H}(h):={\rm Ad}_G(h)
\vert_{\mathfrak q}\,\,(h\in H)$, where ${\rm Ad}_G$ is the adjoint representation of $G$.  
Denote by 
$\mu$ the equivalence class of the complexification of $\rho_{G/H}$.  
Denote by $B_{\mathfrak g}$ the Killing form of 
$\mathfrak g$.  For $\beta\in\triangle'_+$, we set 
$m_{\beta}:={\rm dim}\,\mathfrak p_{\beta},\,\,m^V_{\beta}:=
{\rm dim}(\mathfrak p_{\beta}\cap\mathfrak q)$ and 
$m^H_{\beta}:={\rm dim}(\mathfrak p_{\beta}\cap\mathfrak h)$.  
Also, let $\beta=\sum\limits_{i=1}^rn^{\beta}_i\beta_i,\,\,\,(\beta\in\triangle'_+)$.  
Let $Z_0$ be a point of $\mathfrak b$.  
We consider the following two conditions for $Z_0$:
$$({\rm I})\quad\left\{
\begin{array}{l}
\displaystyle{
\vspace{0.2truecm}\beta(Z_0)\equiv 0,\,\frac{\pi}{6}\,\frac{\pi}{3}\,
\frac{\pi}{2},\,\frac{2\pi}{3}\,\frac{5\pi}{6}\,\,
({\rm mod}\,\pi)\,\,\,\,(\forall\,\beta\in\triangle'_+)\,\,\,\&}\\
\hspace{0.5truecm}\displaystyle{
\sum_{\beta\in{\triangle'}^V_+\,\,{\rm s.t.}\,\,
\beta(Z_0)\equiv\,\frac{\pi}{6}\,\,({\rm mod}\,\pi)}
3n^{\beta}_im^V_{\beta}
+\sum_{\beta\in{\triangle'}^V_+\,\,{\rm s.t.}\,\,
\beta(Z_0)\equiv\,\frac{\pi}{3}\,\,({\rm mod}\,\pi)}
n^{\beta}_im^V_{\beta}}\\
\hspace{0.5truecm}\displaystyle{
+\sum_{\beta\in{\triangle'}^H_+\,\,{\rm s.t.}\,\,
\beta(Z_0)\equiv\frac{2\pi}{3}\,({\rm mod}\,\pi)}3n^{\beta}_im^H_{\beta}
+\sum_{\beta\in{\triangle'}^H_+\,\,{\rm s.t.}\,\,
\beta(Z_0)\equiv\frac{5\pi}{6}\,({\rm mod}\,\pi)}n^{\beta}_im^H_{\beta}}\\
\displaystyle{
=\sum_{\beta\in{\triangle'}^V_+\,\,{\rm s.t.}\,\,
\beta(Z_0)\equiv\,\frac{2\pi}{3}\,\,({\rm mod}\,\pi)}
n^{\beta}_im^V_{\beta}
+\sum_{\beta\in{\triangle'}^V_+\,\,{\rm s.t.}\,\,
\beta(Z_0)\equiv\,\frac{5\pi}{6}\,\,({\rm mod}\,\pi)}
3n^{\beta}_im^V_{\beta}}\\
\hspace{0.5truecm}\displaystyle{
+\sum_{\beta\in{\triangle'}^H_+\,\,{\rm s.t.}\,\,
\beta(Z_0)\equiv\frac{\pi}{6}\,({\rm mod}\,\pi)}n^{\beta}_im^H_{\beta}
+\sum_{\beta\in{\triangle'}^H_+\,\,{\rm s.t.}\,\,
\beta(Z_0)\equiv\frac{\pi}{3}\,({\rm mod}\,\pi)}3n^{\beta}_im^H_{\beta}}\\
\hspace{8.8truecm}\displaystyle{(i=1,\cdots,r).}
\end{array}\right.$$
and
$$({\rm II})\quad\left\{
\begin{array}{l}
\displaystyle{
\beta(Z_0)\equiv 0,\,\frac{\pi}{4},\,\frac{\pi}{2},\,\frac{3\pi}{4}\,\,
({\rm mod}\,\pi)\,\,\,\,(\forall\,\beta\in\triangle'_+)\,\,\,\&}\\
\hspace{0.5truecm}\displaystyle{
\sum_{\beta\in{\triangle'}^V_+\,\,{\rm s.t.}\,\,
\beta(Z_0)\equiv\frac{\pi}{4}\,({\rm mod}\,\pi)}n^{\beta}_im^V_{\beta}
+\sum_{\beta\in{\triangle'}^H_+\,\,{\rm s.t.}\,\,
\beta(Z_0)\equiv\frac{3\pi}{4}\,({\rm mod}\,\pi)}n^{\beta}_im^H_{\beta}}\\
\displaystyle{
=\sum_{\beta\in{\triangle'}^V_+\,\,{\rm s.t.}\,\,
\beta(Z_0)\equiv\frac{3\pi}{4}\,({\rm mod}\,\pi)}n^{\beta}_im^V_{\beta}
+\sum_{\beta\in{\triangle'}^H_+\,\,{\rm s.t.}\,\,
\beta(Z_0)\equiv\frac{\pi}{4}\,({\rm mod}\,\pi)}n^{\beta}_im^H_{\beta}}\\
\hspace{8.8truecm}\displaystyle{(i=1,\cdots,r).}
\end{array}\right.$$
Denote by $H_{Z_0}$ the isotropy group of the $H$-action at ${\rm Exp}\,Z_0$.  
For simplicity, we set $L:=H_{Z_0}$ and denote the identity component of $L$ by $L_0$.  
Set $M:=H({\rm Exp}\,Z_0)(=H/L)$ and $\widehat M:=H/L_0$, and define a covering map 
$\psi:\widehat M\to M$ by $\psi(hL_0)=hL\,\,\,(h\in H)$.  
Denote by $\iota$ the inclusion map of $M$ into $G/K$ and set 
$\widehat{\iota}:=\iota\circ\psi$.  
In the sequel, we regard $\widehat M$ as a submanifold in $G/K$ embedded by 
$\widehat{\iota}$.  
Also, denote by $\mathfrak h_{Z_0}$ (or $\mathfrak l$) the Lie algebra of $L$.  
We showed that $M$ is minimal and that 
$\mathfrak h$ admits a natural reductive decomposition 
$\mathfrak h=\mathfrak l+\mathfrak m_{\mathfrak h}$ (see Theorem A in [Koi3] or 
the proof of Theorem A of this paper).  
Furthermore, we ([Koi3]) showed that the induced metric on the submanifold $M$ in $G/K$ 
coincides with the $H$-invariant metric arising from the restriction 
$cB_{\mathfrak g}\vert_{\mathfrak m_{\mathfrak h}\times\mathfrak m_{\mathfrak h}}$ of some 
constant-multiple $cB_{\mathfrak g}$ of $B_{\mathfrak g}$ to 
$\mathfrak m_{\mathfrak h}\times\mathfrak m_{\mathfrak h}$ if one of the following 
conditions holds:

\vspace{0.2truecm}


$({\rm I}_1)\quad$ $({\rm I})$ holds, 
${\triangle'}^V_+\cap{\triangle'}^H_+=\emptyset$, 
$\beta(Z_0)\equiv 0,\,\frac{\pi}{3},\,\frac{2\pi}{3}\,\,({\rm mod}\,\pi)$ for all 
$\beta\in{\triangle'}^V_+$ 

\hspace{1.2truecm}and $\beta(Z_0)\equiv\frac{\pi}{6},\,\frac{\pi}{2},\,\frac{5\pi}{6}\,\,
({\rm mod}\,\pi)$ for all $\beta\in{\triangle'}^H_+$, 

\vspace{0.2truecm}

$({\rm I}_2)\quad$ $({\rm I})$ holds, ${\triangle'}^V_+\cap{\triangle'}^H_+=\emptyset$, 
$\beta(Z_0)\equiv 0,\,\frac{\pi}{6},\,\frac{5\pi}{6}\,\,({\rm mod}\,\pi)$ for all 
$\beta\in{\triangle'}^V_+$

\hspace{1.2truecm}and $\beta(Z_0)\equiv\frac{\pi}{3},\,\frac{\pi}{2},\,\frac{2\pi}{3}\,\,
({\rm mod}\,\pi)$ for all $\beta\in{\triangle'}^H_+$, 

\newpage


$({\rm II}_1)\quad$ $({\rm II})$ holds, ${\triangle'}^V_+\cap{\triangle'}^H_+=\emptyset$, 
$\beta(Z_0)\equiv 0,\,\frac{\pi}{4},\,\frac{3\pi}{4}\,\,({\rm mod}\,\pi)$ for all 
$\beta\in{\triangle'}^V_+$ 

\hspace{1.2truecm}and $\beta(Z_0)\equiv\frac{\pi}{4},\,\frac{\pi}{2},\,\frac{3\pi}{4}\,\,
({\rm mod}\,\pi)$ for all $\beta\in{\triangle'}^H_+$

\noindent
(see Theorems $C\sim F$ in [Koi3]).  
Here we note that, when $G$ is simple, 
there exists an inner automorphism $\rho$ of $G$ with $\rho(K)=H$ 
by Proposition 4.39 of [I2].  
Denote by $H^s$ the semi-simple part of $H$ and 
$\mathfrak h^s$ the Lie algebra of $H^s$.  
Let $k$ be the positive integer defined by 
$$k:=\left\{
\begin{array}{ll}
1 & (G/H:{\rm Hermite}\,\,{\rm type})\\
3 & (G/H:{\rm quarternionic}\,\,{\rm Kaehler}\,\,{\rm type})\\
0 & (G/H:{\rm other}).
\end{array}\right.\leqno{(1.1)}$$
Easily we can show $H=S^k\cdot H^s$, where $k$ is as above.  
Denote by $H^s_{Z_0}$ the isotropy group of $H^s$ at ${\rm Exp}\,Z_0$.  
For simplicity, we set $L^s:=H^s_{Z_0}$ and denote the identity component of $L^s$ by 
$L^s_0$.  Denote by $\mathfrak l^s$ the Lie algebra of $L^s$ and 
$\mathfrak z$ the center of $\mathfrak h$ and 
$\mathfrak z_{\mathfrak h}(\mathfrak b)$ the centralizer of $\mathfrak b$ in 
$\mathfrak h$.  
In the case where $G/H$ is of Hermite type or quarternionic Kaehler type, 
we assume that ${\rm cohom}\,H={\rm rank}\,G/K$ holds, where ${\rm cohom}\,H$ is 
the cohomogeneity of the $H$-action.  
From this assumption, $\mathfrak b$ is a maximal abelian subspace of $\mathfrak p$ and 
hence 
$\mathfrak z_{\mathfrak h}(\mathfrak b)=\mathfrak z_{\mathfrak k\cap\mathfrak h}
(\mathfrak b)$.  Also, we have $\mathfrak z\subset\mathfrak z_{\mathfrak h}(\mathfrak b)$ 
(see Page 92 of [Tak]).  
Hence we obatain 
$\mathfrak z\subset\mathfrak z_{\mathfrak k\cap\mathfrak h}(\mathfrak b)$.  
On the other hand, according to $(3.1)$, we have 
$\mathfrak z_{\mathfrak k\cap\mathfrak h}(\mathfrak b)\subset\mathfrak l$.  
Therefore, we obtain $\mathfrak z\subset\mathfrak l$ and hence $L=S^k\cdot L^s$.  
From this relation, it follows that $M=H/L=H^s/L^s$ and that $\widehat M=H/L_0=H^s/L_0^s$.  
Define a covering map $\psi:{\widehat M}\to M$ by 
$\psi(hL^s_0)=hL^s\,\,\,(h\in H^s)$.  Denote by $\iota$ the inclusion map of $M$ into 
$G/K$ and set ${\widehat{\iota}}:=\iota\circ\psi$.  
In the sequel, we regard $\widehat M$ as a submanifold in $G/K$ embedded by 
$\widehat{\iota}$.  Clearly we have 
$\mathfrak h^s=\mathfrak l^s+\mathfrak m_{\mathfrak h}$.  
Let $(\rho_{H^s}^S)_{Z_0}:L^s\to GL(T^{\perp}_{{\rm Exp}\,Z_0}M^s)$ the slice 
representation of the $H^s$-action at ${\rm Exp}\,Z_0$, 
where $T^{\perp}_{{\rm Exp}\,Z_0}M$ is the normal space of $M$ at ${\rm Exp}\,Z_0$.  
Set $\mathfrak m:=(\exp\,Z_0)_{\ast}^{-1}(T_{{\rm Exp}\,Z_0}M)$ and 
$\mathfrak m^{\perp}:=(\exp\,Z_0)_{\ast}^{-1}(T^{\perp}_{{\rm Exp}\,Z_0}M)$.  
Let $I(\exp\,Z_0):G\to G$ be the inner automorphism by $\exp\,Z_0$.  
Easily we can show $I(\exp(-Z_0))(L^s)\subset K$ and hence 
$${\rm Ad}_G(\exp(-Z_0))(\mathfrak l^s)\subset\mathfrak k.\leqno{(1.2)}$$
Also we can show 
$${\rm Ad}_G(\exp\,Z_0)(\mathfrak m^{\perp})\subset\mathfrak q.
\leqno{(1.3)}$$
See $(3.5)$ about the proof of $(1.3)$.  
Set $\mathfrak q^s:=\mathfrak z+\mathfrak q$.  
Also, let $\rho_{G/H}:H\to GL(\mathfrak q)$ be the isotropy representation of 
$G/H$ and $\rho_{G/H^s}:H^s\to GL(\mathfrak q^s)$ the isotropy representation of 
$G/H^s$.  Define the representation 
$\sigma_{Z_0}:L^s_0\to GL({\rm Ad}_G(\exp\,Z_0)(\mathfrak m^{\perp}))$ by 
$$\sigma_{Z_0}({\it l})(w):=(\rho_{G/H^s}({\it l}))(w)\qquad
({\it l}\in L^s_0,\,w\in{\rm Ad}_G(\exp\,Z_0)(\mathfrak m^{\perp})).$$
Under the identification of $T^{\perp}_{{\rm Exp}\,Z_0}M$ and 
${\rm Ad}_G(\exp\,Z_0)(\mathfrak m^{\perp})$, 
the restriction $(\rho_H^S)_{Z_0}\vert_{L_0^s}$ of $(\rho_H^S)_{Z_0}$ to $L_0^s$ 
is identified with $\sigma_{Z_0}$.  
We regard $({\rm Ad}_G(\exp\,Z_0)(\mathfrak m^{\perp}))^{\bf c}$ 
as a $L_0^s$-module associated with the complexification 
$\sigma_{Z_0}^{\bf c}:L_0^s\to 
GL(({\rm Ad}_G(\exp\,Z_0)(\mathfrak m^{\perp}))^{\bf c})$ of $\sigma_{Z_0}$.  
Denote by $\mu$ the equivalence class of the complexification 
$\rho^{\bf c}_{G/H}:H\to GL(\mathfrak q^{\bf c})$ of $\rho_{G/H}$ and 
$\mu\vert_{H^s}$ the equivalence class of of the restriction 
$\rho^{\bf c}_{G/H}\vert_{H^s}$ of $\rho^{\bf c}_{G/H}$ to $H^s$.  

In this paper, we prove the following result.  

\vspace{0.5truecm}

\noindent
{\bf Theorem A.} {\sl Let $G/K$ be an irreducible simply connected symmetric 
space of compact type, 
$H\curvearrowright G/K$ a Hermann action and $Z_0$ an element of $\mathfrak b$ such that 
$(H,Z_0)$ satisfies one of the above conditions $({\rm I}_1),\,({\rm I}_2)$ or 
$({\rm II}_1)$.  
Furthermore, assume that 
${\rm cohom}\,H={\rm rank}\,G/K$ holds.  
Let $M,\,\,\widehat M,\,\,H^s$ and $L_0$ 
be the quantities defined for $(H,Z_0)$ as above.  
Then the orbit $M$ (hence $\widehat M$) is minimal (but not totally geodesic) and the index 
$i(\widehat M)$ of $\widehat M$ is given by 
$$i(\widehat M)=\sum_{\lambda\in D_{G/H}}m_{\lambda}\cdot
{\rm dim}\,{\rm Hom}_{L_0^s}(V_{\rho_{\lambda}},({\rm Ad}_G(\exp\,Z_0)
(\mathfrak m^{\perp}))^{\bf c}).$$
Here $D_{G/H}:=\{\lambda\in D(H^s)\,\vert\,a_{\lambda}\,>\,a_{\mu\vert_{H^s}}\}$, where 
$a_{\lambda}$ (resp. $a_{\mu\vert_{H^s}}$) is the eigenvalue of the Casimir operator of 
an irreducible complex representation belonging to 
$\lambda$ (resp. $\mu\vert_{H^s}$) with respect to 
$B_{\mathfrak g}\vert_{\mathfrak h^s\times\mathfrak h^s}$, 
$V_{\rho_{\lambda}}$ is the representation space of an irreducible 
representation $\rho_{\lambda}$ belonging to $\lambda$, $m_{\lambda}$ is 
the dimension of $V_{\rho_{\lambda}}$ and 
${\rm Hom}_{L_0^s}(V_{\rho_{\lambda}},
({\rm Ad}_G(\exp\,Z_0)({\mathfrak m}^{\perp}))^{\bf c})$ is the 
$L_0^s$-module of all $L_0^s$-homomorphisms from $V_{\rho_{\lambda}}$ to 
$({\rm Ad}_G(\exp\,Z_0)({\mathfrak m}^{\perp}))^{\bf c}$.
}

\vspace{0.5truecm}

\noindent
{\it Remark 1.1.} 
(i) In general, we have $i(M)\leq i(\widehat M)$.  In particular, if $L$ is connected, then 
we have $M=\widehat M$.  

(ii) If $G/H$ is of Hermite-type, then the isotrpy representation $\rho_{G/H}$ of 
$G/H$ is an irreducible complex representation of $H$ and, 
when its equivalence class is denoted by $\nu$, we have $\mu=\nu\oplus\nu$ 
and $a_{\mu}=a_{\nu}$.  

\vspace{0.5truecm}

In the final section, we give examples of a Hermann action 
$H\curvearrowright G/K$ and $Z_0\in\mathfrak b$ as in Theorem A and 
calculate the indices of $\widehat M$ for some of the examples 
by using Theorem A.  

\section{Basic notions and facts} 
In this section, we recall some basic notions and facts.  

\newpage


\noindent
{\bf Jacobi operators}

\vspace{0.1truecm}

Let $f:(M,g)\hookrightarrow(\widetilde M,\widetilde g)$ be a minimal isometric 
immersion of a compact Riemannian manifold $(M,g)$ into another Riemannian 
manifold $(\widetilde M,\widetilde g)$.  Denote by $T^{\perp}M$ the normal 
bundle of $f$ and $\Gamma(T^{\perp}M)$ the space of all normal vector 
fields of $f$.  Also, denote by $\nabla$(resp. $\nabla^{\perp}$) the 
Levi-Civita connection of $g$ (resp. the normal connection of $f$) and 
$A$ the shape tensor of $f$.  
Let $f_t$ ($-\varepsilon<t<\varepsilon$) be a $C^{\infty}$-family of 
immersions of $M$ into $\widetilde M$ with $f_0=f$, where $\varepsilon$ is 
a positive number.  Define a map $F:M\times(-\varepsilon,\varepsilon)\to
\widetilde M$ by $F(x,t):=f_t(x)$ ($(x,t)\in 
M\times(-\varepsilon,\varepsilon)$).  Denote by 
${\rm Vol}(M,f_t^{\ast}\widetilde g)$ the volume of 
$(M,f_t^{\ast}\widetilde g)$ and $dv$ the volume element of $g$, where 
$f_t^{\ast}\widetilde g$ is the metric induced form $\widetilde g$ by $f_t$.  
Then we have the following second variational formula:
$$\left.\frac{d^2}{dt^2}\right\vert_{t=0}{\rm Vol}(M,f_t^{\ast}\widetilde g)
=\int_M\widetilde g\left({\cal J}\left(
F_{\ast}\left(\left.\frac{\partial}{\partial t}\right\vert_{t=0}\right)_{\perp}
\right),
F_{\ast}\left(\left.\frac{\partial}{\partial t}\right\vert_{t=0}\right)_{\perp}
\right)dv$$
(see Theorem 3.2.2 in [S]).  
Here $F_{\ast}$ is the differential of $F$, 
$(\cdot)_{\perp}$ is the normal component of $(\cdot)$, ${\cal J}$ is 
the Jacobi operator of $f$ (or $M$), which is defined by 
${\cal J}:=-\triangle^{\perp}+{\cal R}-{\cal A}\,
(:\Gamma(T^{\perp}M)\to\Gamma(T^{\perp}M))$ (where $\triangle^{\perp}$ is 
the rough Laplacian operator defined by $\nabla$ and $\nabla^{\perp}$, 
${\cal A}$ is defined by $g({\cal A}(v),w)={\rm Tr}(A_v\circ A_w)$ 
($v,w\in\Gamma(T^{\perp}M)$) and ${\cal R}$ is defined by 
$g({\cal R}(v),w)=-{\rm Tr}(R(\cdot,v)w)$ ($v,w\in\Gamma(T^{\perp}M)$)).  
Set $E_{\lambda}^{\perp}:=\{v\in\Gamma(T^{\perp}M)\,\vert\,{\cal J}(v)
=\lambda v\}$ for each $\lambda\in{\Bbb R}$.  The dimension of 
$\sum\limits_{\lambda<0}E_{\lambda}^{\perp}$ (resp. $E_0^{\perp}$) is called 
the {\it index} (resp. {\it nullity}) of $f$ (or $M$).  

\vspace{0.3truecm}

\noindent
{\bf The eigenvalues of the Casimir operators}

\vspace{0.1truecm}

For a compact Lie group $H$, denote by $D(H)$ the set of all 
equivalence classes of (finite dimensional) irreducible complex 
representations of $H$.  
Fix an ${\rm Ad}(H)$-invariant inner product $\langle\,\,,\,\,\rangle$ of the Lie algebra 
$\mathfrak h$ of $H$.  Let $\rho$ be an irreducible complex representation of $H$.  
The Casimir operator $C_{\rho}$ of $\rho$ with respect to $\langle\,\,,\,\,\rangle$ 
is defined by 
$C_{\rho}:=\sum\limits_{i=1}^m\rho_{\ast e}(e_i)^2$, 
where $(e_1,\cdots,e_m)$ is an orthonormal base of $\mathfrak h$ with respect to 
$\langle\,\,,\,\,\rangle$ and $e$ is the identity element of $H$.  
Assume that $H$ is semi-simple and connected.  
Fix a Cartan subalgebra $\widetilde{\mathfrak a}$ of the Lie algebra 
$\mathfrak h$ of $H$.  Let $\triangle$ be the root system of $\mathfrak h$ 
with respect to $\widetilde{\mathfrak a}$, $\triangle_+$ the positive root 
system of $\triangle$ under some lexicographic ordering of the dual space 
${\widetilde{\mathfrak a}}^{\ast}$ of $\widetilde{\mathfrak a}$ and 
$\Pi=\{\alpha_1,\cdots,\alpha_r\}$ be a simple root system of $\triangle_+$.  
Define $\Lambda_i\in{\widetilde{\mathfrak a}}^{\ast}$ ($i=1,\cdots,r$) by 
$\displaystyle{\frac{2\langle\alpha_j,\Lambda_i\rangle}{\langle\alpha_j,
\alpha_j\rangle}=\delta_{ij}}$ ($1\leq i,j\leq r$).  
It is known that an injection of $D(H)$ into ${\Bbb Z}_+\{\Lambda_1,\cdots,\Lambda_r\}\,
(:=\{\sum_{i=1}^rz_i\Lambda_i\,\vert\,z_i\in{\Bbb Z}_+\})$ is given by assigning the highest 
weight of $\rho$ to each $\lambda=[\rho]\in D(H)$, where $[\rho]$ is the equivalence 
class of an irreducible complex representation $\rho$ of $H$.  
Denote by $\widehat{D(H)}$ the image of this injection.  Then 
The quotient group ${\Bbb Z}_+\{\Lambda_1,\cdots,\Lambda_r\}/\widehat{D(H)}$ is isomorphic 
to the fundamental group $\pi_1(H)$ of $H$.  
Denote by $(z_1,\cdots,z_r)$ the equivalence class of the irreducible complex representation 
of $H$ corresponding to $\sum\limits_{i=1}^rz_i\Lambda_i$.  
If $H$ is simple, then we have 
$C_{\rho}=a_{\rho}{\rm id}_{\mathfrak h}$ for some $a_{\rho}\in{\Bbb R}$ 
(by Schur's lemma), where ${\rm id}_{\mathfrak h}$ is the identity transformation of 
$\mathfrak h$.  According to the Freudenthal's formula, we have 
$$a_{\rho}=-\langle\Lambda,\Lambda+\sum_{\alpha\in\triangle_+}\alpha\rangle,
\leqno{(2.1)}$$
where $\Lambda$ is the highest weight of $\rho$.  

\vspace{0.3truecm}

\noindent
{\bf Irreducible complex representations of $T^r,\,Spin(2r)$ and $Spin(2r+1)$}

\vspace{0.1truecm}

For each $(m_1,\cdots,m_r)\in{\Bbb Z}^r$, an ireducible complex representation 
$\rho$ of $r$-dimensional torus group $T^r(=SO(2)^r=U(1)^r)$ is defined by 
$$\rho(z_1,\cdots,z_r)(w):=z_1^{m_1}\cdots z_r^{m_r}w\,\,\,
((z_1,\cdots,z_r)\in T^r=U(1)^r,\,\,w\in{\Bbb C}).$$
Denote by $(m_1-\cdots-m_r)$ the equivalence class of this representation.  
Let $D(T^r)$ be the set of all the equivalence classes of irreducible complex 
representations of $T^r$.  Then it is known that 
$D(T^r)=\{(z_1-\cdots -z_r)\,\vert\,(z_1,\cdots, z_r)\in{\Bbb Z}^r\}$ holds 
(see [KO] for example).  

Let $H=Spin(2r)$ or $Spin(2r+1)$, and ${\widetilde{\mathfrak a}}^{\ast}$ and 
$\Pi=\{\alpha_1,\cdots,\alpha_r\}$ be as above.  
Also, let $\{\beta_1,\cdots,\beta_r\}$ be the base of 
${\widetilde{\mathfrak a}}^{\ast}$ defined by 
$\alpha_i=\beta_i-\beta_{i+1}\,\,(i=1,\cdots,r-1)$ and 
$$\left\{\begin{array}{ll}
\alpha_r=\beta_{r-1}+\beta_r & (H=Spin(2r))\\
\alpha_r=\beta_r & (H=Spin(2r+1)).
\end{array}\right.$$
For an irrecducible complex representation $\rho$ of $H$, the highest weight 
$\Lambda$ of $\rho$ is expressed as $\Lambda=\sum\limits_{i=1}^rm_i\beta_i$ for some 
$(m_1,\cdots,m_r)\in{\Bbb Z}^r+\{(0,\cdots,0),(\frac 12,\cdots,\frac 12)\}$.  
Then we denote the equivalence class of $\rho$ 
by $(m_1\,\cdots\,m_r)^{\bullet}$.  
It is known that 
$$D(H)=\{(m_1\,\cdots\,m_r)^{\bullet}\,\vert\,
(m_1,\cdots,m_r)\in{\Bbb Z}^r+\{(0,\cdots,0),(\frac 12,\cdots,\frac 12)\}\}$$
and that 
$$D(H/\{\pm 1\})=\{(m_1\,\cdots\,m_r)^{\bullet}\,\vert\,
(m_1,\cdots,m_r)\in{\Bbb Z}^r\},$$
where $H/\{\pm 1\}=SO(2r)$ or $SO(2r+1)$ (see Chapter 9 of [KO] for example).  

\vspace{0.3truecm}

\noindent
{\bf The canonical connection}

\vspace{0.1truecm}

Let $H/L$ be a reductive homogeneous space and $\mathfrak h=\mathfrak l
+\mathfrak m$ be a reductive decomposition (i.e., 
$[\mathfrak l,\mathfrak m]\subset\mathfrak m$), where $\mathfrak h$ 
(resp. $\mathfrak l$) is the Lie algebra of $H$ (resp. $L$).  
Also, let $\pi:P\to H/L$ be a principal $G$-bundle, where $G$ is a Lie group.  
Assume that $H$ acts on $P$ as $\pi(h\cdot u)=h\cdot\pi(u)$ for any 
$u\in P$ and any $h\in H$.  
Then there uniquely exists a connection $\omega$ of 
$P$ such that, for any $X\in\mathfrak m$ and any $u\in P$, 
$t\mapsto(\exp\,tX)(u)$ is a horizontal curve with respect to 
$\omega$, where $\exp$ is the exponential map of $H$.  This connection 
$\omega$ is called the {\it canonical connection} of $P$ associated with 
the reductive decomposition $\mathfrak h=\mathfrak l+\mathfrak m$.  

\vspace{0.3truecm}

\noindent
{\bf The rough Laplacian operator with respect to the canonical connection}

\vspace{0.1truecm}

Let $H$ be a Lie group and $H/L$ be a reductive homogeneous space 
with a reductive decomposition $\mathfrak h=\mathfrak l+\mathfrak m$, 
where $\mathfrak h$ is the Lie algebra of $H$.  
The subspace $\mathfrak m$ is identified with $T_{eL}(H/L)$.  Let $B$ be an 
${\rm Ad}(H)$-invariant inner product of $\mathfrak h$ such that 
$\mathfrak h=\mathfrak l+\mathfrak m$ is an orthogonal decomposition with 
respect to $B$.  
Denote by $\langle\,\,,\,\,\rangle$ the $H$-invariant metric on $H/L$ induced 
from $B\vert_{\mathfrak m\times\mathfrak m}$ and $\nabla$ the Levi-Civita 
connection of $\langle\,\,,\,\,\rangle$.  
Let $\pi:H\to H/L$ be the natural projection, 
$\sigma:L\to GL(W)$ a unitary representation of $L$ and 
$E_{\sigma}:=H\times_{\sigma(L)}W$ the associated complex vector bundle of 
the $L$-bundle $\pi:H\to H/L$ with respect to $\sigma$.  The Lie group $H$ 
acts on $H$ and $H/L$ naturally.  Also, each $h(\in H)$ gives a linear 
isomorphism of $W$ onto the fibre $(E_{\sigma})_{\pi(h)}$.  
Denote by $\Gamma(E_{\sigma})$ the space of all sections of $E_{\sigma}$ 
and set $C^{\infty}(H,W)_{\sigma}:=\{f\in C^{\infty}(H,W)\,\vert\,
f(h{\it l})=\sigma({\it l}^{-1})f(h)\,(\forall\,h\in H,\forall\,{\it l}\in L)
\}$, where $C^{\infty}(H,W)$ is the space of all $W$-valued 
$C^{\infty}$-functions on $H$.  Define a map 
$\Psi:\Gamma(E_{\sigma})\to C^{\infty}(H,W)_{\sigma}$ by 
$\Psi(\xi)(h)=h^{-1}\cdot\xi_{\pi(h)}$ ($\xi\in\Gamma(E_{\sigma}),\,h\in H$).  
This map $\Psi$ is a linear isomorphism preserving the $H$-action.  
Take an orthonormal base $(e_1,\cdots,e_m)$ of $\mathfrak h$ with respect to 
$B$ with $e_i\in\mathfrak l$ ($i=1,\cdots,n$) and $e_b\in\mathfrak m$ 
($b=n+1,\cdots,m$), where $n:={\rm dim}\,\mathfrak l$.  Let 
${\cal C}_H\,(:C^{\infty}(H,W)\to C^{\infty}(H,W))$ be the Casimir 
differential operator of $H$ with respect to $B$, that is, 
${\cal C}_H(f)=\sum\limits_{i=1}^m\widetilde e_i(\widetilde e_if)$, where 
$\widetilde e_i$ is the left-invariant vector field induced from $e_i$.  
Also, let ${\cal C}_{\sigma}$ be the Casimir operator of $\sigma$ with 
respect to $B\vert_{\mathfrak l\times\mathfrak l}$.  For 
$f\in C^{\infty}(H,W)_{\sigma}$, we can show 
${\cal C}_H(f)={\cal C}_{\sigma}\circ f+\sum\limits_{b=n+1}^m
\widetilde e_b(\widetilde e_bf)$.  Let $\nabla^{\omega}$ be the connection of 
$E_{\sigma}$ induced from the canonical connection $\omega$ of 
$\pi:H\to H/L$ with respect to the reductive decomposition 
$\mathfrak h=\mathfrak l+\mathfrak m$ and $\triangle^{E_{\sigma}}$ the 
rough Laplacian operator of $E_{\sigma}$ with respect to $\nabla^{\omega}$ 
and $\nabla$.  Set $\widetilde{\triangle^{E_{\sigma}}}:=\Psi\circ
\triangle^{E_{\sigma}}\circ\Psi^{-1}$.  Then we have 
$\widetilde{\triangle^{E_{\sigma}}}f=\sum\limits_{b=n+1}^m
\widetilde e_b(\widetilde e_bf)$ ($f\in C^{\infty}(H,W)_{\sigma}$) 
by Proposition 2.3 of [O].  Furthermore, by Corollary 2.5 of [O], 
we have the following relation.  

\vspace{0.5truecm}

\noindent
{\bf Lemma 2.1([O]).} {\sl For each $f\in C^{\infty}(H,W)_{\sigma}$, we have 
$$\widetilde{\triangle^{E_{\sigma}}}f
={\cal C}_H(f)-{\cal C}_{\sigma}\circ f.$$}

\section{Proof of Theorem A}
In this section, we shall prove Theorem A.  
We use the notations in Introduction.  
Let $(H,Z_0)$ be as in the statement of Theorem A.  
Denote by $\langle\,\,,\,\,\rangle$ the $G$-invariant metric of $G/K$ 
induced from $B_{\mathfrak g}\vert_{\mathfrak p\times\mathfrak p}$.  
We shall described some subspaces stated in Introduction 
explicitly.  
Set 
$${\triangle'}^V_{Z_0}:=\{\beta\in{\triangle'}^V_+\,\vert\,\beta(Z_0)\equiv 0\,\,
({\rm mod}\,\pi)\}$$
and 
$${\triangle'}^H_{Z_0}:=\{\beta\in{\triangle'}^H_+\,\vert\,\beta(Z_0)\equiv \frac{\pi}{2}
\,\,({\rm mod}\,\pi)\}.$$
Clearly the Lie algebra $\mathfrak l$ is given by 
$$\mathfrak l=
\mathfrak z_{\mathfrak k\cap\mathfrak h}(\mathfrak b)
+\sum_{\beta\in{\triangle'}^V_{Z_0}}
(\mathfrak k_{\beta}\cap\mathfrak h)
+\sum_{\beta\in{\triangle'}^H_{Z_0}}(\mathfrak p_{\beta}\cap\mathfrak h).
\leqno{(3.1)}$$
Easily we can show that $\mathfrak m_{\mathfrak h}$ and 
$\mathfrak m_{\mathfrak h}^s$ are given by 
$$\mathfrak m_{\mathfrak h}=
\mathfrak z_{\mathfrak p\cap\mathfrak h}(\mathfrak b)
+\sum_{\beta\in{\triangle'}^V_+\setminus{\triangle'}^V_{Z_0}}
(\mathfrak k_{\beta}\cap\mathfrak h)
+\sum_{\beta\in{\triangle'}^H_+\setminus{\triangle'}^H_{Z_0}}
(\mathfrak p_{\beta}\cap\mathfrak h).
\leqno{(3.2)}$$
From these relations, it follws that the decompositions 
$\mathfrak h=\mathfrak l+\mathfrak m_{\mathfrak h}$ and 
$\mathfrak h^s=\mathfrak l^s+\mathfrak m_{\mathfrak h}$ are reductive, respectively.  
Easily we can show that $\mathfrak m$ is given by 
$$\mathfrak m=
\mathfrak z_{\mathfrak p\cap\mathfrak h}(\mathfrak b)+
\sum_{\beta\in{\triangle'}^V_+\setminus{\triangle'}^V_{Z_0}}
(\mathfrak p_{\beta}\cap\mathfrak q)
+\sum_{\beta\in{\triangle'}^H_+\setminus{\triangle'}^H_{Z_0}}
(\mathfrak p_{\beta}\cap\mathfrak h)
\leqno{(3.3)}$$
and hence 
$$\mathfrak m^{\perp}=\mathfrak b
+\sum_{\beta\in{\triangle'}^V_{Z_0}}(\mathfrak p_{\beta}\cap\mathfrak q)
+\sum_{\beta\in{\triangle'}^H_{Z_0}}(\mathfrak p_{\beta}\cap\mathfrak h).
\leqno{(3.4)}$$
Furthermore, we can show 
$$\begin{array}{l}
\hspace{0.5truecm}\displaystyle{{\rm Ad}_G(\exp\,Z_0)(\mathfrak m^{\perp})}\\
\displaystyle{=\mathfrak b
+\sum_{\beta\in{\triangle'}^V_{Z_0}}
\left(\cos({\rm ad}(Z_0))(\mathfrak p_{\beta}\cap\mathfrak q)
+\sin({\rm ad}(Z_0))(\mathfrak p_{\beta}\cap\mathfrak q)\right)}\\
\hspace{0.5truecm}\displaystyle{+\sum_{\beta\in{\triangle'}^H_{Z_0}}
\left(\cos({\rm ad}(Z_0))(\mathfrak p_{\beta}\cap\mathfrak h)
+\sin({\rm ad}(Z_0))(\mathfrak p_{\beta}\cap\mathfrak h)\right)}\\
\displaystyle{=\mathfrak b
+\sum_{\beta\in{\triangle'}^V_{Z_0}}(\mathfrak p_{\beta}\cap\mathfrak q)
+\sum_{\beta\in{\triangle'}^H_{Z_0}}(\mathfrak k_{\beta}\cap\mathfrak q)
\qquad(\subset\mathfrak q).}
\end{array}
\leqno{(3.5)}$$

\vspace{0.5truecm}

\noindent
{\it Proof of Theorem A.} 
From $(3.1)$ and $(3.2)$, we have 
$B_{\mathfrak g}(\mathfrak l,\mathfrak m_{\mathfrak h})=0$.  
By imitating the discussion in the proof of Theorem A in [Koi3], we can show that 
$(\psi^{\ast}g_I)_{eL}=cB_{\mathfrak g}\vert_{\mathfrak m_{\mathfrak h}\times
\mathfrak m_{\mathfrak h}}$, where 
$c=\frac34$ in case of $({\rm I}_1)$, $c=\frac14$ in case of $({\rm I}_2)$ 
and $c=\frac12$ in case of $({\rm II}_1)$.  
Let $\omega$ be the canonical connection of the principal $L_0$-bundle 
$\pi:H^s\to H^s/L^s_0(={\widehat M})$ with respect to the reductive 
decomposition $\mathfrak h^s=\mathfrak l^s+\mathfrak m_{\mathfrak h}$ and 
$F^{\perp}({\widehat M})$ the normal frame bundle of ${\widehat M}$.  
Note that $F^{\perp}({\widehat M})$ is identified with the induced bundle 
$\psi^{\ast}(F^{\perp}(M))\,\,(\subset{\widehat M}\times F^{\perp}(M))$ of $F^{\perp}(M)$ 
by $\psi$.  Define a map $\eta:H^s\to F^{\perp}({\widehat M})$ by 
$\eta(h)=(hL_0^s,h_{\ast}u_0)$ ($h\in H^s$), where $u_0$ is a fixed normal frame of 
$M$ at ${\rm Exp}\,Z_0$.  This map $\eta$ is an embedding.  
By identifying $H^s$ with $\eta(H^s)$, we regard $\pi:H^s\to H^s/L^s_0(={\widehat M})$ 
as a subbundle of $F^{\perp}({\widehat M})$.  
Denote by the same symbol $\omega$ the connection of $F^{\perp}({\widehat M})$ 
induced from $\omega$ and $\nabla^{\omega}$ the linear connection on 
$T^{\perp}{\widehat M}$ associated with $\omega$.  
Denote by $\nabla^{\perp}$ the normal connection of the submanifold ${\widehat M}$.  
By imitating the discussion in the proof of Theorem A in [Koi3], we can show that 
$\nabla^{\omega}=\nabla^{\perp}$.  
Denote by 
$E_{\sigma_{Z_0}}$ the associated vector bundle 
$H^s\times_{\sigma_{Z_0}}{\rm Ad}_G(\exp\,Z_0)(\mathfrak m^{\perp})$ 
of the principal $L_0$-bundle $\pi:H^s\to H^s/L_0$ with respect to $\sigma_{Z_0}$, 
where $\sigma_{Z_0}$ is as stated in Introduction.  
Since $\sigma_{Z_0}$ is identified with $(\rho_H^S)_{Z_0}\vert_{L_0^s}$ as stated in 
Introduction, $E_{\sigma_{Z_0}}$ is identified with the normal bundle 
$T^{\perp}{\widehat M}$ of ${\widehat M}$ under 
the correspendence $h\cdot v\,\,\leftrightarrow\,\,
(hL_0,h_{\ast}((\exp\,Z_0)_{\ast}({\rm Ad}_G(\exp(-Z_0))(v))))$ 
($h\in H^s,\,\,v\in{\rm Ad}_G(\exp\,Z_0)(\mathfrak m^{\perp})$).  
Also we note that $T^{\perp}{\widehat M}$ is identified with the induced bundle 
$\psi^{\ast}(T^{\perp}M)\,\,(\subset{\widehat M}\times T^{\perp}M)$ of 
$T^{\perp}M$ by $\psi^s$.  
Let $\Psi:\Gamma(E_{\sigma_{Z_0}})\to C^{\infty}(H^s,
{\rm Ad}_G(\exp\,Z_0)(\mathfrak m^{\perp}))_{\sigma_{Z_0}}$ be a 
diffeomorphism defined in the previous section.  
Denote by $\nabla$ the Levi-Civita connection of $\psi^{\ast}g_I$.  
Since $\phi^{\ast}g_I$ coincides with the $H^s$-invariant metric induced from 
$cB_{\mathfrak g}\vert_{\mathfrak m_{\mathfrak h}\times\mathfrak m_{\mathfrak h}}$ and 
$\nabla^{\omega}=\nabla^{\perp}$, 
it follows from Lemma 2.1 that the rough Laplacian operator 
$\triangle^{\perp}$ of $E_{\sigma_{Z_0}}$ with respect to $\nabla^{\perp}$ and 
$\nabla$ satisfies 
$$\begin{array}{r}
(\Psi\circ\triangle^{\perp}\circ\Psi^{-1})(f)
={\cal C}_{H^s}(f)-{\cal C}_{\sigma_{Z_0}}\circ f\\
(f\in C^{\infty}(H^s,{\rm Ad}_G(\exp\,Z_0)(\mathfrak m^{\perp}))_{\sigma_{Z_0}},
\end{array}
\leqno{(3.6)}$$
where ${\cal C}_{H^s}$ is the Casimir differential operator of $H^s$ 
with respect to $cB_{\mathfrak g}\vert_{\mathfrak h^s\times\mathfrak h^s}$ 
and ${\cal C}_{\sigma_{Z_0}}$ is the Casimir operator of 
$\sigma_{Z_0}$ with respect to $cB_{\mathfrak g}\vert_{\mathfrak l\times\mathfrak l}$.  
Let ${\cal R}$ and ${\cal A}$ be the operators defined for ${\widehat M}$ 
in similar to ${\cal R}$ and ${\cal A}$ stated in the previous section, respectively.  
Then, by using Lemma 4.1 of [I1], we can show 
$$\begin{array}{r}
\displaystyle{
(\Psi\circ{\cal R}\circ\Psi^{-1})(f)=\sum_{i=1}^n[
(e_i)_{\mathfrak p},
[
(e_i)_{\mathfrak p},f]]_{{\rm Ad}_G(\exp\,Z_0)(\mathfrak m^{\perp})}}\\
\displaystyle{(f\in C^{\infty}(H^s,{\rm Ad}_G(\exp\,Z_0)(\mathfrak m^{\perp}))
_{\sigma_{Z_0}}))}
\end{array}
\leqno{(3.7)}$$
and 
$$
\begin{array}{r}
\displaystyle{
(\Psi\circ{\cal A}\circ\Psi^{-1})(f)
=-\sum_{i=1}^n[
(e_i)_{\mathfrak k},
[
(e_i)_{\mathfrak k},f]]_{{\rm Ad}_G(\exp\,Z_0)(\mathfrak m^{\perp})}}\\
\displaystyle{
(f\in C^{\infty}(H^s,{\rm Ad}_G(\exp\,Z_0)(\mathfrak m^{\perp}))_{\sigma_{Z_0}}))}
\end{array}
\leqno{(3.8)}$$
where $(e_1,\cdots,e_n)$ is an orthonormal base of $\mathfrak m_{\mathfrak h}$ 
with repect to $cB_{\mathfrak g}\vert_{\mathfrak m_{\mathfrak h}\times
\mathfrak m_{\mathfrak h}}$, and 
$(\cdot)_{\mathfrak k},\,\,(\cdot)_{\mathfrak p}$ and 
$(\cdot)_{{\rm Ad}_G(\exp\,Z_0)(\mathfrak m^{\perp})}$ is the 
$\mathfrak k$-component, $\mathfrak p$-component and \newline
${\rm Ad}_G(\exp\,Z_0)(\mathfrak m^{\perp})$-component of $(\cdot)$, respectively.  
From $(3.6),\,(3.7)$ and $(3.8)$, the Jacobi operator ${\cal J}$ of ${\widehat M}$ 
is given by 
$$\begin{array}{r}
\displaystyle{(\Psi\circ{\cal J}\circ\Psi^{-1})(f)
=-{\cal C}_{H^s}(f)+
{\cal C}_{\rho_{G/H^s}}\circ f}\\
\displaystyle{
(f\in C^{\infty}(H^s,{\rm Ad}_G(\exp\,Z_0)(\mathfrak m^{\perp}))_{\sigma_{Z_0}}).}
\end{array}
\leqno{(3.9)}$$
Easily we can show 
$$\rho_{G/H^s}(h)={\rm id}_{\mathfrak z}\oplus\rho_{G/H}(h)$$
for any $h\in H^s$, and hence 
$${\cal C}_{\rho_{G/H^s}}=0_{\mathfrak z}\oplus
\frac{a_{\mu\vert_{H^s}}}{c}{\rm id}_{\mathfrak q},$$
where $0_{\mathfrak z}$ is the zero map from 
$\mathfrak z$ to one-self and $a_{\mu\vert_{H^s}}$ is as in the statement of Theorem A.  
Hence we have 
$$\begin{array}{c}
\displaystyle{(\Psi\circ{\cal J}\circ\Psi^{-1})(f)
=-{\cal C}_{H^s}(f)+\frac{a_{\mu\vert_{H^s}}}{c}f}\\
\displaystyle{(f\in C^{\infty}(H^s,{\rm Ad}_G(\exp\,Z_0)(\mathfrak m^{\perp}))
_{\sigma_{Z_0}}).}
\end{array}
\leqno{(3.10)}$$
Let $\lambda(=[\rho_{\lambda}])$ be an element of $D(H^s)$.  
Define a map $\eta_{\rho_{\lambda}}:V_{\rho_{\lambda}}\otimes
{\rm Hom}_{L_0^s}(V_{\rho_{\lambda}},$\newline
$({\rm Ad}_G(\exp\,Z_0)(\mathfrak m^{\perp}))^{\bf c})\to\,
C^{\infty}(H^s,({\rm Ad}_G(\exp\,Z_0)(\mathfrak m^{\perp}))^{\bf c})
_{(\sigma_{Z_0})^{\bf c}}$ 
by 
$$\begin{array}{c}
(\eta_{\rho_{\lambda}}(v\otimes\phi))(h):=\phi(\rho_{\lambda}(h^{-1})(v))\\
(v\in V_{\rho_{\lambda}},\,\,\phi\in{\rm Hom}_{L_0^s}(V_{\rho_{\lambda}},
({\rm Ad}_G(\exp\,Z_0)(\mathfrak m^{\perp}))^{\bf c}),\,\,h\in H^s).
\end{array}$$
This map $\eta_{\rho_{\lambda}}$ is injective.  
Denote by $E_{(\sigma_{Z_0})^{\bf c}}$ the associated complex vector bundle 
$H^s\times_{(\sigma_{Z_0})^{\bf c}}({\rm Ad}_G(\exp\,Z_0)
(\mathfrak m^{\perp}))^{\bf c}$ 
of $\pi:H^s\to H^s/L_0^s$ with respect to $(\sigma_{Z_0})^{\bf c}$, which is identified 
with the complexification $(T^{\perp}{\widehat M})^{\bf c}$ of $T^{\perp}{\widehat M}$.  
Define a diffeomorphism $(\Psi^s)^{\bf c}:\Gamma(E_{(\sigma_{Z_0})^{\bf c}})\to 
C^{\infty}(H^s,({\rm Ad}_G(\exp\,Z_0)(\mathfrak m^{\perp}))^{\bf c})
_{(\sigma_{Z_0})^{\bf c}}$ by 
$\Psi^{\bf c}(\xi)(h):=h^{-1}\cdot\xi_{\pi(h)}$ 
($\xi\in\Gamma(E_{(\sigma_{Z_0})^{\bf c}}),\,h\in H^s)$.  
Set $\Gamma_{\lambda}((T^{\perp}{\widehat M})^{\bf c}):=(\Psi^{\bf c})^{-1}
(\eta_{\rho_{\lambda}}(V_{\rho_{\lambda}}\otimes{\rm Hom}_{L_0^s}(V_{\rho_{\lambda}},$
\newline
$({\rm Ad}_G(\exp\,Z_0)(\mathfrak m^{\perp}))^{\bf c}))$.  
Then, according to Peter-Weyl theorem for vector bundles (see Page P173 of 
[B]), $\sum_{\lambda\in D(H^s)}\Gamma_{\lambda}((T^{\perp}{\widehat M})^{\bf c})$ 
(direct sum) is uniformly dense in $\Gamma((T^{\perp}{\widehat M})^{\bf c})$ with respect 
to the uniformly topology.  Also, it follows from $(3.10)$ that 
$$({\cal J})^{\bf c}(f)=\frac{a_{\mu\vert_{H^s}}-a_{\lambda}}{c}f\quad
(f\in\Gamma_{\lambda}((T^{\perp}{\widehat M})^{\bf c})).
\leqno{(3.11)}$$
From this relation, we have 
$$i({\widehat M})=\sum_{\lambda\in D_{G/H}}m_{\lambda}\cdot
{\rm dim}\,{\rm Hom}_{L_0^s}(V_{\rho_{\lambda}},({\rm Ad}_G(\exp\,Z_0)
({\mathfrak m}^{\perp}))^{\bf c}),$$
where $D_{G/H}$ is as in the statement of Theorem A.  
This completes the proof.  
\hspace{0.5truecm}q.e.d.

\vspace{0.5truecm}

\section{Examples} 
In this section, we give examples of a Hermann action 
$H\curvearrowright G/K$ and $Z_0\in\mathfrak b$ as in Theorem A and 
calculate the index of the minimal orbit $M:=H({\rm Exp}\,Z_0)$ for some of 
the examples by using Theorem A.  
We use the notations in Introduction.  
First we give examples of $(H,Z_0)$ satisfying the condition $({\rm I}_1)$.  

\vspace{0.5truecm}

\noindent
{\it Example 1.} We consider the isotropy action of $SU(3n+3)/SO(3n+3)$.  
Then we have $\triangle=\triangle'$, which is of $({\mathfrak a}_{3n+2})$-type.  
Also, we have $\triangle'_+={\triangle'}^V_+$ and hence 
${\triangle'}^H_+=\emptyset$.  Let $\Pi=\{\beta_1,\cdots,\beta_{3n+2}\}$ be 
a simple root system of $\triangle'_+$, where we order 
$\beta_1,\cdots,\beta_{3n+2}$ as 
the Dynkin diagram of $\triangle'_+$ is as in Figure 1, 
$\triangle'_+=\{\beta_i+\cdots+\beta_j\,\vert\,1\leq i,j\leq 3n+2\}$.  
For any $\beta\in\triangle'_+$, we have $m_{\beta}=1$.  
Let $Z_0$ be the point of $\mathfrak b$ defined by 
$\beta_{n+1}(Z_0)=\beta_{2n+2}(Z_0)=\frac{\pi}{3}$ and $\beta_i(Z_0)=0$ 
($i\in\{1,\cdots,3n+2\}\setminus\{n+1,2n+2\}$).  
This point $Z_0$ satisfies the condition $({\rm I}_1)$ (see Section 4 of [Koi3]).  

\vspace{0.5truecm}

\centerline{
\unitlength 0.1in
\begin{picture}( 11.7000,  1.8000)(  3.3000, -7.2000)
%
\special{pn 8}%
\special{ar 440 590 50 50  0.0000000 6.2831853}%
%
\special{pn 8}%
\special{ar 650 590 50 50  0.0000000 6.2831853}%
%
\special{pn 8}%
\special{ar 1450 590 50 50  0.0000000 6.2831853}%
%
\special{pn 8}%
\special{pa 500 590}%
\special{pa 600 590}%
\special{fp}%
%
\special{pn 8}%
\special{pa 700 590}%
\special{pa 760 590}%
\special{fp}%
%
\special{pn 8}%
\special{pa 1330 590}%
\special{pa 1400 590}%
\special{fp}%
%
\special{pn 8}%
\special{pa 810 590}%
\special{pa 1240 590}%
\special{dt 0.045}%
\put(3.3000,-7.0000){\makebox(0,0)[lt]{{\scriptsize$\beta_1$}}}%
\put(5.7000,-7.0000){\makebox(0,0)[lt]{{\scriptsize$\beta_2$}}}%
\put(13.6000,-7.2000){\makebox(0,0)[lt]{{\scriptsize$\beta_{3n+2}$}}}%
\end{picture}%
\hspace{0truecm}}

\vspace{0.5truecm}

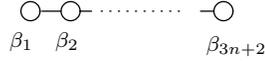
\centerline{{\bf Figure 1.}}

\vspace{0.5truecm}

\noindent
{\it Example 2.} 
We consider the isotropy action of $SU(6n+6)/Sp(3n+3)$.  
Then we have $\triangle=\triangle'$, which is of $({\mathfrak a}_{3n+2})$-type.  
Also, we have $\triangle'_+={\triangle'}^V_+$ and hence 
${\triangle'}^H_+=\emptyset$.  
Let $\Pi=\{\beta_1,\cdots,\beta_{3n+2}\}$ be 
a simple root system of $\triangle'_+$, where we order 
$\beta_1,\cdots,\beta_{3n+2}$ as above.  We have $m_{\beta}=4$ for any 
$\beta\in\triangle'_+$.  
Let $Z_0$ be the point of $\mathfrak b$ defined by 
$\beta_{n+1}(Z_0)=\beta_{2n+2}(Z_0)=\frac{\pi}{3}$ and $\beta_i(Z_0)=0$ 
($i\in\{1,\cdots,3n+2\}\setminus\{n+1,2n+2\}$).  
This point $Z_0$ satisfies the condition $({\rm I}_1)$ (see Section 4 of [Koi3]).  

\vspace{0.5truecm}

\noindent
{\it Example 3.} 
We consider the isotropy action of $SU(3)/S(U(1)\times U(2))$ 
($2$-dimensional complex projective space).  
Then we have $\triangle=\triangle'$, which is of $({\mathfrak bc}_1)$-type.  
Also, we have $\triangle'_+={\triangle'}^V_+$ and hence 
${\triangle'}^H_+=\emptyset$.  
Let $\Pi=\{\beta\}$ be a simple root system of $\triangle'_+$.  
We have $\triangle'_+=\{\beta,2\beta\}$ and $m_{\beta}=2$ and $m_{2\beta}=4$.  
Let $Z_0$ be the point of $\mathfrak b$ defined by 
$\beta(Z_0)=\frac{\pi}{3}$.  This point $Z_0$ satisfies the condition 
$({\rm I}_1)$ (see Section 4 of [Koi3]).  

\vspace{0.5truecm}

\noindent
{\it Example 4.} 
We consider the isotropy action of $Sp(3n+2)/U(3n+2)$.  
Then we have $\triangle=\triangle'$, which is of $({\mathfrak c}_{3n+2})$-type.  
Also, we have $\triangle'_+={\triangle'}^V_+$ and hence 
${\triangle'}^H_+=\emptyset$.  
Let $\Pi=\{\beta_1,\cdots,\beta_{3n+2}\}$ be 
a simple root system of $\triangle'_+$, where we order 
$\beta_1,\cdots,\beta_{3n+2}$ as the Dynkin diagram of $\triangle'_+$ is as in Figure 2.  
We have $m_{\beta}=1$ for any 
$\beta\in\triangle'_+$.  
Let $Z_0$ be the point of $\mathfrak b$ defined by 
$\beta_{n+1}(Z_0)=\beta_{3n+2}(Z_0)=\frac{\pi}{3}$ and $\beta_i(Z_0)=0$ 
($i\in\{1,\cdots,3n+2\}\setminus\{n+1,3n+2\}$).  
This point $Z_0$ satisfies the condition $({\rm I}_1)$ (see Section 4 of [Koi3]).  

\vspace{0.5truecm}

\centerline{
\unitlength 0.1in
\begin{picture}( 14.0000,  1.6000)( 10.2000, -8.4000)
%
\special{pn 8}%
\special{ar 1140 730 50 50  0.0000000 6.2831853}%
%
\special{pn 8}%
\special{ar 1350 730 50 50  0.0000000 6.2831853}%
%
\special{pn 8}%
\special{ar 2150 730 50 50  0.0000000 6.2831853}%
%
\special{pn 8}%
\special{pa 1200 730}%
\special{pa 1300 730}%
\special{fp}%
%
\special{pn 8}%
\special{pa 1400 730}%
\special{pa 1460 730}%
\special{fp}%
%
\special{pn 8}%
\special{pa 2030 730}%
\special{pa 2100 730}%
\special{fp}%
%
\special{pn 8}%
\special{pa 1510 730}%
\special{pa 1940 730}%
\special{dt 0.045}%
%
\special{pn 8}%
\special{ar 2370 730 50 50  0.0000000 6.2831853}%
%
\special{pn 8}%
\special{pa 2320 710}%
\special{pa 2230 710}%
\special{fp}%
%
\special{pn 8}%
\special{pa 2200 730}%
\special{pa 2270 680}%
\special{fp}%
%
\special{pn 8}%
\special{pa 2230 780}%
\special{pa 2230 780}%
\special{fp}%
%
\special{pn 8}%
\special{pa 2260 790}%
\special{pa 2190 730}%
\special{fp}%
\put(10.2000,-8.2000){\makebox(0,0)[lt]{{\scriptsize$\beta_1$}}}%
\put(12.8000,-8.3000){\makebox(0,0)[lt]{{\scriptsize$\beta_2$}}}%
\put(23.1000,-8.4000){\makebox(0,0)[lt]{{\scriptsize$\beta_{3n+2}$}}}%
\put(19.1000,-8.3000){\makebox(0,0)[lt]{{\scriptsize$\beta_{3n+1}$}}}%
%
\special{pn 8}%
\special{pa 2240 770}%
\special{pa 2340 770}%
\special{fp}%
\end{picture}%
}

\vspace{0.5truecm}

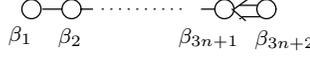
\centerline{{\bf Figure 2.}}

\vspace{0.5truecm}

\noindent
{\it Example 5.} We consider the dual action 
$\rho_1(SO(3))\curvearrowright SU(3)/SO(3)$ of the Hermann action 
$SO_0(1,2)\curvearrowright SL(3,{\Bbb R})/SO(3)$, where 
$\rho_1$ is an inner automorphism of $SU(3)$.  
Then $\triangle=\triangle'$ is of $({\mathfrak a}_2)$-type.  
Let $\Pi=\{\beta_1,\beta_2\}$ be a simple root system of $\triangle'_+$.  
Then we have 
${\triangle'}^V_+=\{\beta_1\}$, ${\triangle'}^H_+=\{\beta_2,\beta_1+\beta_2\}$ 
and hence ${\triangle'}^V_+\cap{\triangle'}^H_+=\emptyset$.  
Also we have $m_{\beta_1}=m_{\beta_2}=m_{\beta_1+\beta_2}=1$.  
Let $Z_0$ be the point of $\mathfrak b$ satisfying 
$(\beta_1(Z_0),\beta_2(Z_0))=(\frac{\pi}{3},-\frac{\pi}{6})$.  
This point $Z_0$ satisfies the condition $({\rm I}_1)$ (see Section 4 of [Koi3]).  

\vspace{0.5truecm}

\noindent
{\it Example 6.} We consider the dual action 
$\rho_2(Sp(3))\curvearrowright SU(6)/Sp(3)$ of the Hermann action 
$Sp(1,2)\curvearrowright SU^{\ast}(6)/Sp(3)$, where 
$\rho_2$ is an inner automorphism of $SU(6)$.  
Then $\triangle=\triangle'$ is of $({\mathfrak a}_2)$-type.  
Let $\Pi=\{\beta_1,\beta_2\}$ be a simple root system of $\triangle'_+$.  
Then we have 
${\triangle'}^V_+=\{\beta_1\}$, ${\triangle'}^H_+=\{\beta_2,\beta_1+\beta_2\}$ 
and hence ${\triangle'}^V_+\cap{\triangle'}^H_+=\emptyset$.  
Also we have $m_{\beta_1}=m_{\beta_2}=m_{\beta_1+\beta_2}=4$.  
Let $Z_0$ be the point of $\mathfrak b$ satisfying 
$(\beta_1(Z_0),\beta_2(Z_0))=(\frac{\pi}{3},-\frac{\pi}{6})$.  
This point $Z_0$ satisfies the condition $({\rm I}_1)$ (see Section 4 of [Koi3]).  

\vspace{0.5truecm}

\noindent
{\it Example 7.} We consider the dual action 
$\rho_3(SU(2))\curvearrowright Sp(2)/U(2)$ of the Hermann action 
$U(1,1)\curvearrowright Sp(2,{\Bbb R})/U(2)$, where 
$\rho_3$ is an inner automorphism of $Sp(2)$.  
Then $\triangle=\triangle'$ is of $({\mathfrak c}_2)$-type.  
Let $\Pi=\{\beta_1,\beta_2\}$ be a simple root system of $\triangle'_+$, where 
we we order $\beta_1,\beta_2$ as 
the Dynkin diagram of $\triangle'_+$ is as in Figure 3.  
Then we have 
${\triangle'}^V_+=\{\beta_2,2\beta_1+\beta_2\}$, 
${\triangle'}^H_+=\{\beta_1,\beta_1+\beta_2\}$ 
and hence ${\triangle'}^V_+\cap{\triangle'}^H_+=\emptyset$.  
Also we have $m_{\beta_1}=m_{\beta_2}=m_{\beta_1+\beta_2}=m_{2\beta_1+\beta_2}=1$.  
Let $Z_0$ be the point of $\mathfrak b$ satisfying 
$(\beta_1(Z_0),\beta_2(Z_0))=(-\frac{\pi}{6},\frac{\pi}{3})$.  
This point $Z_0$ satisfies the condition $({\rm I}_1)$ (see Section 4 of [Koi3]).  

\vspace{0.5truecm}

\centerline{
\unitlength 0.1in
\begin{picture}( 22.3000,  1.5000)(  2.9000, -8.4000)
%
\special{pn 8}%
\special{ar 2140 740 50 50  0.0000000 6.2831853}%
%
\special{pn 8}%
\special{ar 2470 750 50 50  0.0000000 6.2831853}%
%
\special{pn 8}%
\special{pa 2190 740}%
\special{pa 2260 690}%
\special{fp}%
%
\special{pn 8}%
\special{pa 2220 790}%
\special{pa 2220 790}%
\special{fp}%
%
\special{pn 8}%
\special{pa 2250 800}%
\special{pa 2180 740}%
\special{fp}%
\put(23.6000,-8.4000){\makebox(0,0)[lt]{{\scriptsize$\beta_2$}}}%
%
\special{pn 8}%
\special{pa 2230 780}%
\special{pa 2430 780}%
\special{fp}%
%
\special{pn 8}%
\special{pa 2230 720}%
\special{pa 2430 720}%
\special{fp}%
\put(22.7000,-8.4000){\makebox(0,0)[rt]{{\scriptsize$\beta_1$}}}%
\end{picture}%
\hspace{5truecm}}

\vspace{0.5truecm}

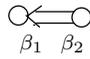
\centerline{{\bf Figure 3.}}


\vspace{0.5truecm}


\vspace{0.5truecm}

\noindent
{\it Example 8.} We consider the dual action 
$\rho_4(Sp(2))\curvearrowright (Sp(2)\times Sp(2))/Sp(2)$ of the Hermann action 
$Sp(1,1)\curvearrowright Sp(2,{\Bbb C})/Sp(2)$, where 
$\rho_4$ is an automorphism of $Sp(2)\times Sp(2)$.  
Then $\triangle=\triangle'$ is of $({\mathfrak c}_2)$-type.  
Let $\Pi=\{\beta_1,\beta_2\}$ be a simple root system of $\triangle'_+$, where 
we we order $\beta_1,\beta_2$ as 
the Dynkin diagram of $\triangle'_+$ is as in Figure 3.  
Then we have 
${\triangle'}^V_+=\{\beta_2,2\beta_1+\beta_2\}$, 
${\triangle'}^H_+=\{\beta_1,\beta_1+\beta_2\}$ 
and hence ${\triangle'}^V_+\cap{\triangle'}^H_+=\emptyset$.  
Also we have $m_{\beta_1}=m_{\beta_2}=m_{\beta_1+\beta_2}=m_{\beta_1+2\beta_2}=2$.  
Let $Z_0$ be the point of $\mathfrak b$ satisfying 
$(\beta_1(Z_0),\beta_2(Z_0))=(-\frac{\pi}{6},\frac{\pi}{3})$.  
This point $Z_0$ satisfies the condition $({\rm I}_1)$ (see Section 4 of [Koi3]).  

\vspace{0.5truecm}

\noindent
{\it Example 9.} We consider the dual action 
$\rho_5(F_4)\curvearrowright E_6/F_4$ of the Hermann action 
$F_4^{-20}\curvearrowright E_6^{-26}/F_4$, where 
$\rho_6$ is an inner automorphism of $E_6$.  
Then $\triangle=\triangle'$ is of $({\mathfrak a}_2)$-type.  
Let $\Pi=\{\beta_1,\beta_2\}$ be a simple root system of $\triangle'_+$.  
Then we have 
${\triangle'}^V_+=\{\beta_1\}$, 
${\triangle'}^H_+=\{\beta_2,\beta_1+\beta_2\}$ 
and hence ${\triangle'}^V_+\cap{\triangle'}^H_+=\emptyset$.  
Also we have $m_{\beta_1}=m_{\beta_2}=m_{\beta_1+\beta_2}=8$.  
Let $Z_0$ be the point of $\mathfrak b$ satisfying 
$(\beta_1(Z_0),\beta_2(Z_0))=(\frac{\pi}{3},-\frac{\pi}{6})$.  
This point $Z_0$ satisfies the condition $({\rm I}_1)$ (see Section 4 of [Koi3]).  

\vspace{0.5truecm}

\noindent
{\it Example 10.} We consider the dual action 
$\rho_6(SO(4))\curvearrowright G_2/SO(4)$ of the Hermann action 
$SL(2,{\Bbb R})\times SL(2,{\Bbb R})\curvearrowright G_2^2/SO(4)$, where 
$\rho_6$ is an inner automorphism of $G_2$.  
Then $\triangle=\triangle'$ is of $({\mathfrak g}_2)$-type.  
Let $\Pi=\{\beta_1,\beta_2\}$ be a simple root system of $\triangle'_+$, where 
we we order $\beta_1,\beta_2$ as 
the Dynkin diagram of $\triangle'_+$ is as in Figure 4.  
Then we have 
${\triangle'}^V_+=\{\beta_1,3\beta_1+2\beta_2\}$, 
${\triangle'}^H_+=\{\beta_2,\beta_1+\beta_2,2\beta_1+\beta_2,3\beta_1+\beta_2\}$ 
and hence ${\triangle'}^V_+\cap{\triangle'}^H_+=\emptyset$.  
Also we have $m_{\beta_1}=m_{\beta_2}=m_{\beta_1+\beta_2}=m_{2\beta_1+\beta_2}
=m_{3\beta_1+\beta_2}=m_{3\beta_1+2\beta_2}=1$.  
Let $Z_0$ be the point of $\mathfrak b$ satisfying 
$(\beta_1(Z_0),\beta_2(Z_0))=(\frac{\pi}{3},-\frac{\pi}{2})$.  
This point $Z_0$ satisfies the condition $({\rm I}_1)$ (see Section 4 of [Koi3]).  

\vspace{0.5truecm}

\noindent
{\it Example 11.} We consider the dual action 
$\rho_7(G_2)\curvearrowright (G_2\times G_2)/G_2$ of the Hermann action 
$G_2^2\curvearrowright G_2^{\bf c}/G_2$, where 
$\rho_7$ is an automorphism of $G_2\times G_2$.  
Then $\triangle=\triangle'$ is of $({\mathfrak g}_2)$-type.  
Let $\Pi=\{\beta_1,\beta_2\}$ be a simple root system of $\triangle'_+$, where 
we we order $\beta_1,\beta_2$ as 
the Dynkin diagram of $\triangle'_+$ is as in Figure 4.  
Then we have 
${\triangle'}^V_+=\{\beta_1,3\beta_1+2\beta_2\}$, 
${\triangle'}^H_+=\{\beta_2,\beta_1+\beta_2,2\beta_1+\beta_2,3\beta_1+\beta_2\}$ 
and hence ${\triangle'}^V_+\cap{\triangle'}^H_+=\emptyset$.  
Also we have $m_{\beta_1}=m_{\beta_2}=m_{\beta_1+\beta_2}=m_{2\beta_1+\beta_2}
=m_{3\beta_1+\beta_2}=m_{3\beta_1+2\beta_2}=2$.  
Let $Z_0$ be the point of $\mathfrak b$ satisfying 
$(\beta_1(Z_0),\beta_2(Z_0))=(\frac{\pi}{3},-\frac{\pi}{2})$.  
This point $Z_0$ satisfies the condition $({\rm I}_1)$ (see Section 4 of [Koi3]).  

\vspace{0.5truecm}

\centerline{
\unitlength 0.1in
\begin{picture}( 22.4000,  1.6000)(  2.9000, -8.4000)
%
\special{pn 8}%
\special{ar 2150 730 50 50  0.0000000 6.2831853}%
%
\special{pn 8}%
\special{ar 2480 740 50 50  0.0000000 6.2831853}%
%
\special{pn 8}%
\special{pa 2200 730}%
\special{pa 2270 680}%
\special{fp}%
%
\special{pn 8}%
\special{pa 2230 780}%
\special{pa 2230 780}%
\special{fp}%
%
\special{pn 8}%
\special{pa 2260 790}%
\special{pa 2190 730}%
\special{fp}%
\put(23.6000,-8.4000){\makebox(0,0)[lt]{{\scriptsize$\beta_2$}}}%
%
\special{pn 8}%
\special{pa 2240 770}%
\special{pa 2440 770}%
\special{fp}%
%
\special{pn 8}%
\special{pa 2240 710}%
\special{pa 2440 710}%
\special{fp}%
\put(22.7000,-8.4000){\makebox(0,0)[rt]{{\scriptsize$\beta_1$}}}%
%
\special{pn 8}%
\special{pa 2210 740}%
\special{pa 2430 740}%
\special{fp}%
\end{picture}%
\hspace{5truecm}}

\vspace{0.5truecm}

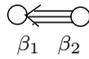
\centerline{{\bf Figure 4.}}

\vspace{0.5truecm}

First we prepare the following lemma.  

\vspace{0.5truecm}

\noindent
{\bf Lemma 4.1.} {\sl Let $G/K,H,L,\theta,\tau$ and $M$ be as in Introduction.  
If both the symmetric space $H/H\cap K$ and the principal orbit of 
the isotropy action of the symmetric space ${\rm Fix}(\theta\circ\tau)_0/H\cap K$ are 
simply connected, then so is also $M$.}

\vspace{0.5truecm}

\noindent
{\it Proof.} Easily we can show $H(eK)=H/H\cap K$ and 
$\exp^{\perp}(T^{\perp}_{eK}H(eK))={\rm Fix}(\theta\circ\tau)_0/H\cap K$, where 
$\exp^{\perp}$ is the normal exponential map of $H(eK)$.  
Let $M'$ be a principal orbit of the isotropy action of 
${\rm Fix}(\theta\circ\tau)_0/H\cap K$.  Then we can show that the focal map of $M'$ onto 
$H(eK)$ is a fibration having a principal orbit of the isotropy action of 
${\rm Fix}(\theta\circ\tau)_0/H\cap K$ as fibre.  
Hence it follows from the assumption that $M'$ is simply connected.  
Let ${\rm pr}$ be the natural projection of $M'$ onto $M$.  
In the case where $M$ is a singular orbit, ${\rm pr}$ is the focal map of $M'$ onto 
$M$ and it is a fibration with connected fibre, where we note that the fibre is the image 
of a principal orbit of the direct sum representation of some 
$s$-representations by the normal exponential map (of $M$) and hence it is connected.  
In the case where $M$ is a principal orbit, ${\rm pr}$ is 
the end-point map (which is a diffeomorphism) of $M'$ onto $M$.  
In both cases, ${\rm pr}$ is a fibration with connected fibre.  
Hence, since $M'$ is simply connected, so is also $M$.  
\hspace{12.1truecm}q.e.d.

\vspace{0.5truecm}

For the representations $\rho_{\lambda_i}$ of $H_i$ ($i=1,\cdots,k$), 
we define the representation $\rho_{\lambda_1}$-$\cdots$-$\rho_{\lambda_k}$ 
of $H_1\times\cdots\times H_k$ by 
$(\rho_{\lambda_1}$-$\cdots$-$\rho_{\lambda_k})(h_1,\cdots,h_k)
(v_1\otimes\cdots\otimes v_k):=\rho_{\lambda_1}(h_1)(v_1)\otimes\cdots\otimes
\rho_{\lambda_k}(h_k)(v_k)$ ($h_i\in H_i,\,\,v_i\in V_{\rho_{\lambda_i}}$) 
(the representation space of $\rho_{\lambda_1}$-$\cdots$-$\rho_{\lambda_k}$ is 
$V_{\rho_{\lambda_1}}\otimes\cdots\otimes V_{\rho_{\lambda_k}}$).  
Denote by $(\lambda_1$-$\cdots$-$\lambda_k)$ the equivalence class of 
$\rho_{\lambda_1}$-$\cdots$-$\rho_{\lambda_k}$.  
By using Theorems A or B, we shall calculate the indices of some of the 
minimal orbits $M=H({\rm Exp}\,Z_0)$ in Examples $1\sim11$.  

First we consider the case of $n=2$ in Example 1 (i.e., the case where 
$G/K=SU(9)/SO(9),\,H=SO(9)$ and $M=SO(9)({\rm Exp}\,Z_0)$ 
($\beta_3(Z_0)=\beta_6(Z_0)=\frac{\pi}{3},\,\beta_i(Z_0)=0\,(i\not=3,6)$)).  
Since $SO(9)$ is simple, we have $H^s=H$.  
The equivalence class $\mu$ of the complexification of 
the isotropy representation of $G/H$ is equal to $(2\,0\,0\,0)$.  
Hence, according to Table 1 in [MP], all of the equivalence classes $\lambda$'s of 
irreducible complex representations of $Spin(9)$ with 
$a_{\lambda}>a_{\mu\vert_{H^s}}$ consist of 
$(0\,0\,0\,0),\,(1\,0\,0\,0),\,(0\,0\,0\,1)$ and $(0\,1\,0\,0)$.  
These equivalence classes $(0\,0\,0\,0),\,(1\,0\,0\,0),\,(0\,0\,0\,1)$ and $(0\,1\,0\,0)$ 
are equal to 
$(0\,0\,0\,0)^{\bullet},\,(1\,0\,0\,0)^{\bullet},\,
(\frac12\,\frac12\,\frac12\,\frac12)^{\bullet}$ and $(1\,1\,0\,0)^{\bullet}$, respectively.  
Hence, $(0\,0\,0\,1)$ is not the equivalence classes of the irreducible complex 
representations of $SO(9)$.  
From this fact and $H=H^s$, we have 
$$D_{G/H}=\{(0\,0\,0\,0),\,(1\,0\,0\,0),\,(0\,1\,0\,0)\}.$$
On the other hand, since $\triangle'=\triangle$ is $(\mathfrak a_8)$-type, 
$\beta_3(Z_0)=\beta_6(Z_0)=\frac{\pi}{3}$ and since $\beta_i(Z_0)=0\,(i\not=3,6)$, 
we have ${\triangle'}^V_{Z_0}=\{\beta_1,\beta_2,\beta_1+\beta_2,\beta_4,\beta_5,
\beta_4+\beta_5,\beta_7,\beta_8,\beta_7+\beta_8\}$, 
${\triangle'}^H_{Z_0}=\emptyset$.  
Also we have $\mathfrak z=\mathfrak z_{\mathfrak k\cap\mathfrak h}(\mathfrak b)=\{0\}$, 
we have 
$\mathfrak l^s=\mathfrak l=\sum\limits_{i\in\{1,\,4,\,7\}}
(\mathfrak h_{\beta_i}+\mathfrak h_{\beta_{i+1}}
+\mathfrak h_{\beta_i+\beta_{i+1}})$.  
Also we have ${\rm dim}\,\mathfrak h_{\beta}=1$ for all $\beta\in{\triangle'}^V_+$.  
Hence we have $\mathfrak l^s=3\mathfrak{so}(3)$.  
%
Hence we have $L_0^s=SO(3)^3$.  Hence, by using Table 2 (the branching rules) in [MP], 
we have the following table:

\vspace{0.3truecm}

$$\begin{tabular}{|c|c|c|}
\hline
{\scriptsize$\lambda$} & {\scriptsize$\lambda\vert_{L_0^s}$} & 
{\scriptsize$m_{\lambda}$}\\
\hline
{\scriptsize$(0\,0\,0\,0)$} & {\scriptsize$(0$-$0$-$0)$} 
& {\scriptsize$1$}\\
\hline
{\scriptsize$(1\,0\,0\,0)$} & {\scriptsize$(2$-$0$-$0)\oplus(0$-$2$-$0)\oplus(0$-$0$
-$2)$} 
& {\scriptsize$9$}\\
\hline
{\scriptsize$(0\,1\,0\,0)$} & {\scriptsize$(2$-$0$-$0)\oplus(0$-$2$-$0)\oplus(0$-$0$-$2)$} 
& {\scriptsize$36$}\\
&{\scriptsize$(2$-$2$-$0)\oplus(2$-$0$-$2)\oplus(0$-$2$-$2)$}&\\
\hline
{\scriptsize$\mu=(2\,0\,0\,0)$} & {\scriptsize$2(0$-$0$-$0)\oplus(4$-$0$-$0)\oplus
(0$-$4$-$0)\oplus(0$-$0$-$4)$} & {\scriptsize$44$}\\
&{\scriptsize$(2$-$2$-$0)\oplus(2$-$0$-$2)\oplus(0$-$2$-$2)$}&\\
\hline
\end{tabular}$$

\centerline{{\bf Table 2.}}

\vspace{0.3truecm}

\noindent
Also we have ${\rm dim}\,\mathfrak m^{\perp}=17$.  
Hence we have 
$$\begin{tabular}{l}
$[(\sigma_{Z_0})^{\bf c}]
=2(0$-$0$-$0)\oplus(4$-$0$-$0)\oplus(0$-$4$-$0)\oplus(0$-$0$-$4)$.\\
\end{tabular}$$
Thus the isomorphicness of the $L_0$-module 
$({\rm Ad}_G(\exp\,Z_0)(\mathfrak m^{\perp}))^{\bf c}$ associated with the representation 
$(\sigma_{Z_0})^{\bf c}$ is analyzed completely.  
Therefore, according to Theorem A, it follows from Table 2 and this fact that 
the index of $\widehat M$ is equal to $2$.  Also, since $Z_0$ belongs to an (open) 
$1$-simplex (which we denote by $\sigma$) of the simplicial complex 
$\overline{\widetilde C}$, $M$ is not stable.  
In fact, when $M$ moves along $\sigma$ as $SO(9)$-orbits, its volume decreases.  
Thus we obtain the following result.  

\vspace{0.5truecm}

\noindent
{\bf Proposition 4.2.} {\sl Let $\Pi=\{\beta_1,\cdots,\beta_8\}$ be the 
simple root system of the positive root system $\triangle_+$ of $SU(9)/SO(9)$ 
($\displaystyle{\mathop{\circ}_{\beta_1}-\mathop{\circ}_{\beta_2}
-\cdots-\mathop{\circ}_{\beta_8}}$) and $Z_0$ the element of $\mathfrak b$ 
with $\beta_3(Z_0)=\beta_6(Z_0)=\frac{\pi}{3}$ and 
$\beta_i(Z_0)=0$ ($i\not=3,6$).  
Then the orbit $M:=SO(9)({\rm Exp}(Z_0))$ of the isotropy action of 
$SU(9)/SO(9)$ is minimal (but not totally geodesic) and 
we have $1\leq i(M)\leq i(\widehat M)=2$, where $\widehat M$ is the above covering 
of $M$.}

\vspace{0.5truecm}

Next we consider the case of $n=1$ in Example 2 (i.e., the case where 
$G/K=SU(12)/Sp(6),\,H=Sp(6),$ and $M=Sp(6)({\rm Exp}\,Z_0)$ 
($\beta_2(Z_0)=\beta_4(Z_0)=\frac{\pi}{3},\,\beta_i(Z_0)=0\,(i\not=2,4)$)).  
Since $Sp(6)$ is simple, we have $H^s=H$.  
The equivalence class $\mu$ of the complexification of 
the isotropy representation of $G/H$ is $(0\,1\,0\,0\,0\,0)$.  
Hence, according to Table 1 in [MP], we have 
$D_{G/H}=\{(0\,0\,0\,0\,0\,0),\,(1\,0\,0\,0\,0\,0)\}$.  
On the other hand, since $\triangle'=\triangle$ is $(\mathfrak a_5)$-type, 
$\beta_2(Z_0)=\beta_4(Z_0)=\frac{\pi}{3}$ and since $\beta_i(Z_0)=0\,(i\not=2,4)$, 
${\triangle'}^V_{Z_0}=\{\beta_1,\beta_3,\beta_5\}$ and 
${\triangle'}^H_{Z_0}=\emptyset$.  
Hence we have 
$\mathfrak l=\mathfrak z_{\mathfrak k\cap\mathfrak h}(\mathfrak b)
+\mathfrak h_{\beta_1}+\mathfrak h_{\beta_3}+\mathfrak h_{\beta_5}$.  
and ${\rm dim}\,\mathfrak h_{\beta_i}=4$ ($i=1,3,5$).  
Also we have $\mathfrak z=\{0\}$ and 
$\mathfrak z_{\mathfrak k\cap\mathfrak h}(\mathfrak b)=6{\mathfrak sp}(1)$ 
From these facts, we have $\mathfrak l^s=3\mathfrak{sp}(2)$.  
Therefore, we have $L_0^s=L_0=Sp(2)^3$.  
Hence, by using Table 2 (the branching rules) in [MP], 
we have the following table:

\vspace{0.2truecm}

$$\begin{tabular}{|c|c|c|}
\hline
{\scriptsize$\lambda$} & {\scriptsize$\lambda\vert_{L^s_0}$} & 
{\scriptsize$m_{\lambda}$}\\
\hline
{\scriptsize$(0\,0\,0\,0\,0\,0)$} & {\scriptsize$(00$-$00$-$00)$} 
& {\scriptsize$1$}\\
\hline
{\scriptsize$(1\,0\,0\,0\,0\,0)$} & {\scriptsize
$(10$-$00$-$00)\oplus(00$-$10$-$00)\oplus(00$-$00$-$10)$} & 
{\scriptsize$12$}\\
\hline
{\scriptsize$\mu=(0\,1\,0\,0\,0\,0)$} & {\scriptsize 
$2(00$-$00$-$00)\oplus(01$-$00$-$00)\oplus(00$-$01$-$00)\oplus(00$-$00$-$01)$} & 
{\scriptsize$65$}\\
& {\scriptsize $\oplus(10$-$10$-$00)\oplus(10$-$00$-$10)
\oplus(00$-$10$-$10)$} & \\
\hline
\end{tabular}$$

\centerline{{\bf Table 3.}}

\vspace{0.2truecm}

\noindent
Also, we have ${\rm dim}\,\mathfrak m^{\perp}=17$.  
Hence we have 
$$\begin{tabular}{l}
$[(\sigma_{Z_0})^{\bf c}]
=2(00$-$00$-$00)\oplus(01$-$00$-$00)\oplus(00$-$01$-$00)\oplus(00$-$00$-$01)$.\\
\end{tabular}$$
Thus the isomorphicness of the $L_0$-module 
$({\rm Ad}_G(\exp\,Z_0)(\mathfrak m^{\perp}))^{\bf c}$ associated with the representation 
$(\sigma_{Z_0})^{\bf c}$ is analyzed completely.  
Therefore, according to Theorem A, it follows from Table 3 and this fact that 
the index of $\widehat M$ is equal to $2$.  
On the other hand, principal orbits of this isotropy action are diffeomorphic to 
$Sp(6)/Sp(1)^6$, which is simply connected.  Also we have $H/H\cap K$ is the one-point set 
because $H=K$.  
Hence, it follows from Lemma 4.1 that $M$ is simply connected, that is, $M=\widehat M$.  
Therefore we obtain the following result.  

\vspace{0.5truecm}

\noindent
{\bf Proposition 4.3.} {\sl Let $\Pi=\{\beta_1,\cdots,\beta_5\}$ be the simple 
root system of the positive root system $\triangle_+$ of $SU(12)/Sp(6)$ 
($\displaystyle{\mathop{\circ}_{\beta_1}-\mathop{\circ}_{\beta_2}
-\mathop{\circ}_{\beta_3}-\mathop{\circ}_{\beta_4}-\mathop{\circ}_{\beta_5}}$) 
and $Z_0$ the element of $\mathfrak b$ with $\beta_2(Z_0)=\beta_4(Z_0)
=\frac{\pi}{3}$ and $\beta_1(Z_0)=\beta_3(Z_0)=\beta_5(Z_0)=0$.  
Then the orbit $M:=Sp(6)({\rm Exp}(Z_0))$ of the isotropy action of 
$SU(12)/Sp(6)$ is minimal (but not totally geodesic) and we have 
$i(M)=2$.}

\vspace{0.5truecm}

Next we consider the case of Example 3 (i.e., 
$G/K=SU(3)/S(U(1)\times U(2)),\,H=S(U(1)\times U(2)),$ and 
$M=S(U(1)\times U(2))({\rm Exp}\,Z_0)$ 
($\beta(Z_0)=\frac{\pi}{3}$)).  Clearly we have $H^s=SU(2)$.  
Since $M$ is a geodesic sphere in $G/K$, it is simply connected and of dimension three.  
Hence we have $L=L_0=U(1)$ and $L_0^s=\{e\}$, where $e$ is the identity element of $G$. 
Since $G/H$ is Hermite-type, the isotropy representation of $G/H$ is regarded as 
an irreducible complex representation of $H\cong U(2)$ and it is equal to 
$(1\,0)^{\bullet}$.  The equivalence class $\mu$ of 
its complexification is equal to $(1\,0)^{\bullet}\oplus(1\,0)^{\bullet}$.  Hence 
we have $\mu\vert_{H^s}=(1\,0)^{\bullet}\vert_{H^s}\oplus(1\,0)^{\bullet}\vert_{H^s}
=(1)\oplus(1)$ and hence $a_{\mu\vert_{H^s}}=a_{(1)}$.  
Hence, according to $(2.3)$ and $(2.18)$ in [MP] and and the Freudenthal's formula, 
we have $D_{G/H}=\{(0)\}$.  
On the other hand, since $\triangle'=\triangle$ is $(\mathfrak{bc}_1)$-type and 
since $\beta(Z_0)=\frac{\pi}{3}$, we have ${\triangle'}^V_{Z_0}={\triangle'}^H_{Z_0}
=\emptyset$.  Also, we have ${\rm dim}(\mathfrak m^{\perp})^{\bf c}=1$.  
According to Theorem A, it follows from these facts and $L^s=L_0^s=\{e\}$ that 
the index of $M$ is equal to $1$.  
Thus we obtain the following result.  

\newpage


\noindent
{\bf Proposition 4.4.} {\sl Let $\triangle_+=\{\beta,2\beta\}$ be the 
positive root system of $SU(3)/S(U(1)\times U(2))$ 
and $Z_0$ the element of $\mathfrak b$ with $\beta(Z_0)=\frac{\pi}{3}$.  
Then the orbit $M:=S(U(1)\times U(2))({\rm Exp}(Z_0))$ of the isotropy action 
of $SU(3)/S(U(1)\times U(2))$ is minimal (but not totally geodesic) 
and we have $i(M)=1$.}

\vspace{0.25truecm}

\noindent
{\it Remark 4.1.} This result has already been proved in [G] in different 
method.  

\vspace{0.25truecm}

Next we consider the case of Example 6 (i.e., 
$G/K=SU(6)/Sp(3),\,H=\rho_2(Sp(3))$ and 
$M=\rho_2(Sp(3))({\rm Exp}\,Z_0)$ 
($\beta_1(Z_0),\beta_2(Z_0))=(\frac{\pi}{3},-\frac{\pi}{6})$)).  
Since $Sp(3)$ is simple, we have $H^s=H=Sp(3)$.  
Since the equivalence class $\mu\vert_{H^s}$ of the complexification of 
the restriction of the isotropy representation of $G/H$ to $H^s$ is $(0\,1\,0)$.  
Hence, according to Table 1 in [MP], we have 
$$D_{G/H}=\{(0\,0\,0),\,\,(1\,0\,0)\}.$$
On the other hand, since $\triangle'=\triangle$ is $(\mathfrak a_2)$-type 
and since $(\beta_1(Z_0),\beta_2(Z_0))=(\frac{\pi}{3},-\frac{\pi}{6})$, 
we have 
${\triangle'}^V_+=\{\beta_1\}$ and that ${\triangle'}^H_+=\{\beta_2,\beta_1+\beta_2\}$.  
we have ${\triangle'}^V_{Z_0}={\triangle'}^H_{Z_0}=\emptyset$.  
Also we have $\mathfrak z_{\mathfrak k}(\mathfrak b)
=\mathfrak z_{\mathfrak k\cap\mathfrak h}(\mathfrak b)=\mathfrak{sp}(1)^3$.  
From these facts, we have $\mathfrak l=\mathfrak{sp}(1)^3$ and hence 
$L_0=L_0^s=Sp(1)^3$.  
Hence, by using Table 2 (the branching rules) in [MP], 
we have the following table:

\vspace{0.2truecm}

$$\begin{tabular}{|c|c|c|}
\hline
{\scriptsize$\lambda$} & {\scriptsize$\lambda\vert_{L^s_0}$} & 
{\scriptsize$m_{\lambda}$}\\
\hline
{\scriptsize$(0\,0\,0)$} & {\scriptsize$(0$-$0$-$0)$} & {\scriptsize$1$}\\
\hline
{\scriptsize$(1\,0\,0)$} & {\scriptsize$(1$-$0$-$0)\oplus(0$-$1$-$0)\oplus(0$-$0$-$1)$} & 
{\scriptsize$6$}\\
\hline
{\scriptsize$\mu\vert_{H^s}=(0\,1\,0)$} & {\scriptsize $2(0$-$0$-$0)\oplus(1$-$1$-$0)
\oplus(1$-$0$-$1)\oplus(0$-$1$-$1)$} & {\scriptsize$14$}\\
\hline
\end{tabular}$$

\centerline{{\bf Table 6.}}

\vspace{0.2truecm}

\noindent
Also, we have ${\rm dim}\,\mathfrak m^{\perp}=2$.  
Hence we have $[(\sigma_{Z_0})^{\bf c}]=2(0$-$0))$.  
Thus the isomorphicness of the $L^s_0$-module 
$({\rm Ad}_G(\exp\,Z_0)(\mathfrak m^{\perp}))^{\bf c}$ associated with the representation 
$(\sigma_{Z_0})^{\bf c}$ is analyzed completely.  
Therefore, according to Theorem A, it follows from Table 6 and this fact that 
the index of $\widehat M$ is equal to $2$.  
On the other hand, we have $H/H\cap K=Sp(3)/Sp(1)\times Sp(2)$ (whcih is simply connected) 
and ${\rm Fix}(\theta\circ\tau)_0/H\cap K=(SU(4)/Sp(2))\times U(1)$.  
The principal orbit of the isotropy action of $(SU(4)/Sp(2))\times U(1)$ is difeomorphic to 
$S^3\times S^3$, which is simply connected.  
Hence, it follows from Lemma 4.1 that $M$ is simply connected, that is, $M=\widehat M$.  
Therefore we obtain the following result.  

\vspace{0.25truecm}

\noindent
{\bf Proposition 4.5.} {\sl Let $\Pi=\{\beta_1,\beta_2\}$ be the simple 
root system of the positive root system $\triangle_+$ of $SU(6)/Sp(3)$ 
($\displaystyle{\mathop{\circ}_{\beta_1}-\mathop{\circ}_{\beta_2}}$) 
and $Z_0$ the element of $\mathfrak b$ with $(\beta_1(Z_0),\beta_2(Z_0))
=(\frac{\pi}{3},-\frac{\pi}{6})$.  
Then the orbit $M:=\rho_2(Sp(3))({\rm Exp}(Z_0))$ of the dual action 
$\rho_2(Sp(3))\curvearrowright SU(6)/Sp(3)$ of 
$Sp(1,2)\curvearrowright SU^{\ast}(6)/Sp(3)$ 
is minimal (but not totally geodesic) and 
we have $i(M)=2$.}

\vspace{0.5truecm}

Next we consider the case of Example 7 (i.e., 
$G/K=Sp(2)/U(2),\,H=\rho_3(U(2))$ and 
$M=\rho_3(U(2))({\rm Exp}\,Z_0)$ 
($\beta_1(Z_0),\beta_2(Z_0))=(\frac{\pi}{3},-\frac{\pi}{6})$)).  
Clearly we have $H^s=SU(2)$.  
Since the equivalence class $\mu\vert_{H^s}$ of the complexification of 
the restriction of the isotropy representation of $G/H$ to $H^s$ is $(2)\oplus(2)$.  
Hence we have $D_{G/H}=\{(0),\,\,(1)\}$, 
On the other hand, since $\triangle'=\triangle$ is $(\mathfrak c_2)$-type 
and since $(\beta_1(Z_0),\beta_2(Z_0))=(-\frac{\pi}{6},\frac{\pi}{3})$, 
we have 
${\triangle'}^V_{Z_0}=\{2beta_1+\beta_2\}$ and ${\triangle'}^H_{Z_0}=\emptyset$.  
Also we have $\mathfrak z_{\mathfrak k}(\mathfrak b)
=\mathfrak z_{\mathfrak k\cap\mathfrak h}(\mathfrak b)=\{0\}$.  
From these facts, we have $\mathfrak l^s=\mathfrak{so}(2)$ and hence 
$L_0^s=SO(2)$.  Denote by $\widetilde{\lambda}$ the trivially extension of 
$\lambda\in D(SU(2))$ to $U(2)$ and $T^2$ a maximal torus of $U(2)$.  
By noticing these facts and using Weyl's character formula 
(see Page 409 of [KO] for example), we have the following table:

\vspace{0.2truecm}

$$\begin{tabular}{|c|c|c|c|c|}
\hline
{\scriptsize$\lambda$} & {\scriptsize$\widetilde{\lambda}$} & 
{\scriptsize$\widetilde{\lambda}\vert_{T^2}$} & {\scriptsize$\lambda\vert_{L_0^s}$} & 
{\scriptsize$m_{\lambda}$}\\
\hline
{\scriptsize$(0)$} & {\scriptsize$(0\,0)^{\bullet}$} & {\scriptsize$(0$-$0)$} 
& {\scriptsize$(0)$} & {\scriptsize$1$}\\
\hline
{\scriptsize$(1)$} & {\scriptsize$(\frac12\,\,(-\frac12))^{\bullet}$} & 
{\scriptsize$(\frac12$-$(-\frac12))\oplus((-\frac12)$-$\frac12)$} 
& {\scriptsize$(\frac12)\oplus(-\frac12)$} & {\scriptsize$2$}\\
\hline
{\scriptsize$(2)$} & {\scriptsize$(1\,\,(-1))^{\bullet}$} & 
{\scriptsize$(1$-$(-1))\oplus(0$-$0)\oplus((-1)$-$1)$} 
& {\scriptsize$(1)\oplus(0)\oplus(-1)$} & {\scriptsize$3$}\\
\hline
\end{tabular}$$

\centerline{{\bf Table 7.}}

\vspace{0.2truecm}

\noindent
Easily we can show that ${\rm dim}\,\mathfrak m^{\perp}=3$ and furthermore 
$[(\sigma_{Z_0})^{\bf c}]=(1)\oplus(0)\oplus(-1)$.  
Therefore, according to Theorem A, it follows from Table 7 and this fact that 
the index of $\widehat M$ is equal to $1$.  Also, since $Z_0$ belongs to an (open) 
$1$-simplex (whcih we denote by $\sigma$) of the simplicial complex 
$\overline{\widetilde C}$, $M$ is not stable.  
In fact, when $M$ moves along $\sigma$ as $\rho_3(U(2))$-orbits, its volume decreases.  
Thus we obtain the following result.  

\vspace{0.3truecm}

\noindent
{\bf Proposition 4.6.} {\sl Let $\Pi=\{\beta_1,\beta_2\}$ be the simple 
root system of the positive root system $\triangle_+$ of $Sp(2)/U(2)$ 
($\displaystyle{\mathop{\circ}_{\beta_1}=>\mathop{\circ}_{\beta_2}}$) 
and $Z_0$ the element of $\mathfrak b$ with $(\beta_1(Z_0),\beta_2(Z_0))
=(-\frac{\pi}{6},\frac{\pi}{3})$.  
Then the orbit $M:=\rho_3(U(2))({\rm Exp}(Z_0))$ of the dual action 
$\rho_3(U(2))\curvearrowright Sp(2)/U(2)$ of 
$U(1,1)\curvearrowright Sp(2,{\Bbb R})/U(2)$ 
is minimal (but not totally geodesic) and we have $i(M)=1$.}

\vspace{0.5truecm}

\vspace{0.5truecm}

Next we consider the case of Example 8 (i.e., 
$G/K=(Sp(2)\times Sp(2))/Sp(2),\,H=\rho_4(Sp(2))$ and 
$M=\rho_4(Sp(2))({\rm Exp}\,Z_0)$ 
($(\beta_1(Z_0),\beta_2(Z_0))=(-\frac{\pi}{6},\frac{\pi}{3})$).  
Clearly we have $H^s=H=Sp(2)$.  
Since the equivalence class $\mu\vert_{H^s}$ of the complexification of 
the restriction of the isotropy representation of $G/H$ to $H^s$ is $(2\,0)$.  
Hence we have $D_{G/H}=\{(0\,0),\,\,(1\,0),\,\,(0\,1)\}$.  
On the other hand, since $\triangle'=\triangle$ is $(\mathfrak b_2)$-type 
and since $(\beta_1(Z_0),\beta_2(Z_0))=(-\frac{\pi}{6},\frac{\pi}{3})$, 
we have 
${\triangle'}^V_{Z_0}=\{\beta_2\}$ and ${\triangle'}^H_{Z_0}=\emptyset$.  
Also we have ${\rm dim}\,\mathfrak z_{\mathfrak k}(\mathfrak b)=
{\rm dim}\,\mathfrak z_{\mathfrak k\cap\mathfrak h}(\mathfrak b)=2$.  
From these facts, we have $\mathfrak l^s=\mathfrak{u}(2)=\mathfrak{so}(2)+\mathfrak{su}(2)$ 
and hence $L_0^s=U(2)$.  
Hence, by using Table 2 (the branching rules) in [MP], 
we have the following table:

\vspace{0.2truecm}

$$\begin{tabular}{|c|c|c|}
\hline
{\scriptsize$\lambda$} & {\scriptsize$\lambda\vert_{L_0^s}$} & 
{\scriptsize$m_{\lambda}$}\\
\hline
{\scriptsize$(0\,0)$} & {\scriptsize$(0$-$0)$} & {\scriptsize$1$}\\
\hline
{\scriptsize$(1\,0)$} & {\scriptsize $(0$-$3)$} & {\scriptsize$4$}\\
\hline
{\scriptsize$(0\,1)$} & {\scriptsize $(0$-$4)$} & {\scriptsize$5$}\\
\hline
{\scriptsize$\mu=\mu\vert_{H^s}=(2\,0)$} & {\scriptsize$(0$-$6)\oplus(0$-$2)$} & 
{\scriptsize$10$}\\
\hline
\end{tabular}$$

\centerline{{\bf Table 8.}}

\vspace{0.2truecm}

\noindent
Also, we have ${\rm dim}\,\mathfrak m^{\perp}=4$.  
Hence we have 
$[(\sigma_{Z_0})^{\bf c}]=(0$-$0)\oplus(0$-$2)$.  
Therefore, according to Theorem A, it follows from Table 9 and this fact that 
the index of $\widehat M$ is equal to $1$.  
On the other hand, we have $H/H\cap K=Sp(2)/Sp(1)\times Sp(1)$ (whcih is simply connected) 
and ${\rm Fix}(\theta\circ\tau)_0/H\cap K=S^3\times S^3$.  
The principal orbit of the isotropy action of $S^3\times S^3$ is difeomorphic to 
$S^2\times S^2$, which is simply connected.  
Hence, it follows from Lemma 4.1 that $M$ is simply connected, that is, $M=\widehat M$.  
Therefore we obtain the following result.  

\vspace{0.3truecm}

\noindent
{\bf Proposition 4.7.} {\sl Let $\Pi=\{\beta_1,\beta_2\}$ be the simple 
root system of the positive root system $\triangle_+$ of $(Sp(2)\times Sp(2))/Sp(2)$ 
($\displaystyle{\mathop{\circ}_{\beta_1}-\mathop{\circ}_{\beta_2}}$) 
and $Z_0$ the element of $\mathfrak b$ with $(\beta_1(Z_0),\beta_2(Z_0))
=(-\frac{\pi}{6},\frac{\pi}{3})$.  
Then the orbit $M:=\rho_5(Sp(2))({\rm Exp}(Z_0))$ of the dual action 
$\rho_5(Sp(2))\curvearrowright Sp(2)/U(2)$ of 
$Sp(1,1)\curvearrowright Sp(2,{\Bbb C})/Sp(2)$ 
is minimal (but not totally geodesic) and 
we have $i(M)=1$.}

\vspace{0.5truecm}

Next we consider the case of Example 9 (i.e., 
$G/K=E_6/F_4,\,H=\rho_6(F_4)$ and 
$M=\rho_5(F_4)({\rm Exp}\,Z_0)$ 
($\beta_1(Z_0),\beta_2(Z_0))=(\frac{\pi}{3},-\frac{\pi}{6})$)).  
Clearly we have $H^s=H=F_4$.  
Since the equivalence class $\mu$ of the complexification of 
the isotropy representation of $G/H$ is $(0\,0\,0\,1)$, 
we have $D_{G/H}=\{(0\,0\,0\,0)\}$.  
On the other hand, since $\triangle'=\triangle$ is $(\mathfrak a_2)$-type 
and since $(\beta_1(Z_0),\beta_2(Z_0))=(\frac{\pi}{3},-\frac{\pi}{6})$, 
we have 
${\triangle'}^V_{Z_0}={\triangle'}^H_{Z_0}=\emptyset$.  
Also we have $z_{\mathfrak k}(\mathfrak b)=
z_{\mathfrak k\cap\mathfrak h}(\mathfrak b)=\mathfrak{so}(8)$.  
From these facts, we have $\mathfrak l^s=\mathfrak{so}(8)$ 
and hence $L_0^s=SO(8)$.  
Hence, by using Table 2 (the branching rules) in [MP], 
we have the following table:

\vspace{0.2truecm}

$$\begin{tabular}{|c|c|c|}
\hline
{\scriptsize$\lambda$} & {\scriptsize$\lambda\vert_{L_0^s}$} & 
{\scriptsize$m_{\lambda}$}\\
\hline
{\scriptsize$(0\,0\,0\,0)$} & {\scriptsize$(0\,0\,0\,0)$} & {\scriptsize$1$}\\
\hline
{\scriptsize$(0\,0\,0\,1)$} & {\scriptsize $2(0\,0\,0\,0)\oplus(1\,0\,0\,0)\oplus(0\,0\,1\,0)
\oplus(0\,0\,0\,1)$} & {\scriptsize$26$}\\
\hline
\end{tabular}$$

\centerline{{\bf Table 9.}}

\vspace{0.2truecm}

\noindent
Also, we have ${\rm dim}\,\mathfrak m^{\perp}=2$.  
Hence we have 
$[(\sigma_{Z_0})^{\bf c}]=2(0\,0\,0\,0)$.  
Therefore, according to Theorem A, it follows from Table 9 and this fact that 
the index of $\widehat M$ is equal to $2$.  
On the other hand, we have $H/H\cap K=F_4/{\rm Spin}(9)$ (whcih is simply connected) 
and ${\rm Fix}(\theta\circ\tau)_0/H\cap K=S^9\times S^1$.  
The principal orbit of the isotropy action of $S^9\times S^1$ is difeomorphic to 
$S^8$, which is simply connected.  
Hence, it follows from Lemma 4.1 that $M$ is simply connected, that is, $M=\widehat M$.  
Therefore we obtain the following result.  

\vspace{0.5truecm}

\noindent
{\bf Proposition 4.8.} {\sl Let $\Pi=\{\beta_1,\beta_2\}$ be the simple 
root system of the positive root system $\triangle_+$ of $E_6/F_4$ 
($\displaystyle{\mathop{\circ}_{\beta_1}-\mathop{\circ}_{\beta_2}}$) 
and $Z_0$ the element of $\mathfrak b$ with $(\beta_1(Z_0),\beta_2(Z_0))
=(\frac{\pi}{3},-\frac{\pi}{6})$.  
Then the orbit $M:=\rho_6(F_4)({\rm Exp}(Z_0))$ of the dual action 
$\rho_5(F_4)\curvearrowright E_6/F_4$ of 
$F_4^{-20}\curvearrowright E_6^{-26}/F_4$ 
is minimal (but not totally geodesic) and 
we have $i(M)=2$.}

\vspace{1truecm}

\centerline{{\bf References}}

\vspace{0.5truecm}

{\small 
\noindent
[B] R. Bott, The index theorem for homogeneous differential operators, 
Differential and 

Combinatorial Topology, Princeton University Press, 1965, 167-187.  

\noindent
[BCO] J. Berndt, S. Console and C. Olmos, Submanifolds and Holonomy, 
Research Notes 

in Mathematics 434, CHAPMAN $\&$ HALL/CRC Press, Boca Raton, London, New York 

Washington, 2003.

\noindent
[Co] L. Conlon, Remarks on commuting involutions, Proc. Amer. Math. Soc. 
{\bf 22} (1969) 

255-257.

\noindent
[G] T. Gotoh, The nullity of compact minimal real hypersurfaces in a 
complex projective 

space, Tokyo J. Math. {\bf 17} (1994) 201-209.

\noindent
[GT] O. Goertsches and G. Thorbergsson, 
On the Geometry of the orbits of Hermann actions, 

Geom. Dedicata {\bf 129} (2007) 101-118.

\noindent
[He] S. Helgason, 
Differential geometry, Lie groups and Symmetric Spaces, Academic Press, 

New York, 1978.

\noindent
[HPTT] E. Heintze, R.S. Palais, C.L. Terng and G. Thorbergsson, 
Hyperpolar actions on 

symmetric spaces, Geometry, topology and physics for Raoul Bott 
(ed. S. T. Yau), 

Conf. Proc. Lecture Notes Geom. Topology {\bf 4}, Internat. Press, Cambridge, 
MA, 1995 

pp214-245.

\noindent
[HTST] D. Hirohashi, H. Tasaki, H. Song and R. Takagi, Minimal orbits 
of the isotropy 

groups of symmetric spaces of compact type, Differential Goem. Appl. {\bf 13} 
(2000) 

167-177.

\noindent
[I1] O. Ikawa, Equivariant minimal immersions of compact Riemannian 
homogeneous spaces 

into compact Riemannian homogeneous spaces, Tsukuba J. Math. {\bf 17} 
(1993) 169-188.  

\noindent
[I2] O. Ikawa, The geometry of symmetric triad and orbit spaces of Hermann 
actions, 

J. Math. Soc. Japan {\bf 63} (2011) 79-139.  

\noindent
[IST] O. Ikawa, T. Sakai and H. Tasaki, Orbits of Hermann actions, Osaka J. 
Math. {\bf 38} 

(2001) 923-930.  

\noindent
[Ki] T. Kimura, Stability of certain reflective submanifolds in compact 
symmetric spaces, 

Tsukuba J. Math. {\bf 32} (2008) 361-382.

\noindent
[KT] T. Kimura and M.S. Tanaka, Stability of certain minimal submanifolds 
in compact 

symmetric spaces, Differential Geom. Appl. {\bf 27} (2009) 23-33.

\noindent
[KO] T. Kobayashi and T. Oshima, Lie Groups and Representations, 
Iwanami, Tokyo 

2005 (in Japanese).

\noindent
[Koi1] N. Koike, Actions of Hermann type and proper complex equifocal 
submanifolds, 

Osaka J. Math. {\bf 42} (2005) 599-611.

\noindent
[Koi2] N. Koike, Collapse of the mean curvature flow for equifocal 
submanifolds, Asian J. 

Math. {\bf 15} (2011) 101-128.

\noindent
[Koi3] N. Koike, Examples of certain kind of minimal orbits of Hermann actions, 
to appear 

in Hokkaido Math. J.

\noindent
[Kol] A. Kollross, A Classification of hyperpolar and cohomogeneity one 
actions, Trans. 

Amer. Math. Soc. {\bf 354} (2001) 571-612.

\noindent
[MP] W.G. Mckay and J. Patera, Tables of dimensions, indices, and branching 
rules for 

reresentations of simple Lie algebras, Lecture Notes in Pure and 
Applied Mathematics 

vol. 69, Marcel Dekker, Inc., New York and Baesel (1981).  

\noindent
[N] T. Nagura, On the Jacobi differential operators associated to minimal 
isometric immer-

sions of symmetric spaces, I, II, III, 
Osaka J. Math. {\bf 18} (1981) 115-145, 
{\bf 19} (1982) 79-124, 

{\bf 19} (1982) 241-281.  

\noindent
[O] Y. Ohnita, On stability of minimal submanifolds in compact symmetric 
spaces, Compo-

sitio Math. {\bf 64} (1987) 157-189.

\noindent
[PT] R.S. Palais and C.L. Terng, Critical point theory and submanifold 
geometry, Lecture 

Notes in Math. {\bf 1353}, Springer, Berlin, 1988.

\noindent
[S] J. Simons, Minimal varieties in riemannian manifolds, Ann. of Math. {\bf 88} 
(1968) 62-105.

\noindent
[Tak] M. Takeuchi, "Modern theory of spherical functions (in Japanese)", Iwanami, 
Tokyo, 

1975.

\noindent
[Tan] M.S. Tanaka, Stability of minimal submanifolds in symmetric spaces, 
Tsukuba J. 

Math. {\bf 19} (1995) 27-56.

\noindent
[TT] C.L. Terng and G. Thorbergsson, 
Submanifold geometry in symmetric spaces, J. Diff-

erential Geometry {\bf 42} (1995) 665-718.

\noindent
[Y1] I. Yokota, Groups and Topology, Shokabou, Tokyo 1971 (in Japanese).

\noindent
[Y2] I. Yokota, Groups and Representations, Shokabou, Tokyo 1973 (in Japanese).

\vspace{0.5truecm}

{\small 
\rightline{Department of Mathematics, Faculty of Science}
\rightline{Tokyo University of Science, 1-3 Kagurazaka}
\rightline{Shinjuku-ku, Tokyo 162-8601 Japan}
\rightline{(koike@ma.kagu.tus.ac.jp)}
}

\end{document}